\documentclass[onefignum,onetabnum]{siamart220329}

\usepackage{blkarray} 
\usepackage{comment}
\usepackage{lipsum}
\usepackage{amsfonts}
\usepackage{graphicx}
\usepackage{epstopdf}
\usepackage{algorithmic}
\usepackage{framed}
\ifpdf
  \DeclareGraphicsExtensions{.eps,.pdf,.png,.jpg}
\else
  \DeclareGraphicsExtensions{.eps}
\fi
\DeclareMathAlphabet{\mathpzc}{OT1}{pzc}{m}{it}
\newcommand{\myrightleftarrows}[1]{\mathrel{\substack{\xrightarrow{#1} \\[-.9ex] \xleftarrow{#1}}}}

\newsiamremark{remark}{Remark}
\newsiamremark{hypothesis}{Hypothesis}
\crefname{hypothesis}{Hypothesis}{Hypotheses}
\newsiamthm{claim}{Claim}
\newtheorem*{thm*}{Theorem}

\headers{Saddle-node bifurcations in chemical reaction networks}{N. Vassena}

\title{Structural conditions for saddle-node bifurcations in chemical reaction networks\thanks{Submitted to the editors on \today.
\funding{This work was funded by the Deutsche Forschung Gesellschaft (German Research Foundation). Project n.512355535}}}

\author{Nicola Vassena\thanks{Universit\"{a}t Leipzig,
  \email{nicola.vassena@uni-leipzig.de}.}}

\usepackage{amsopn}

\ifpdf
\hypersetup{
  pdftitle={Saddle-node bifurcations in chemical reaction networks},
  pdfauthor={N. Vassena}
}
\fi

\begin{document}

\maketitle

% REQUIRED
\begin{abstract}
Motivated by investigating multistationarity in biochemical systems, we address saddle-node bifurcations for chemical reaction networks endowed with general kinetics. At positive equilibria, we identify structural network conditions that guarantee the bifurcation behavior and we develop a method to identify the proper bifurcation parameters. As a relevant example, we explicitly provide such bifurcation parameters for Michaelis-Menten and Hill kinetics. Examples of application include reversible feedback cycles, the central carbon metabolism of \emph{Escherichia coli}, and autocatalytic networks.
\end{abstract}

% REQUIRED

% REQUIRED

\section{Introduction}

Multistationarity is the property of a chemical system to exhibit two or more distinct equilibria, under identical conditions, and it has been proposed as an explanation for many epigenetic processes, including cell differentiation: see the groundbreaking work \cite{ThomKauf01} by Thomas and Kaufman, and the many biological references therein. Hence, it is no surprise that investigating multistationarity for chemical systems has become a hot topic. See among others the works by Soul{\'e} \cite{Soule2003}, Craciun and Feinberg \cite{Crafei05, Crafei06}, Mincheva and Roussel \cite{MR07}, Banaji and Craciun \cite{BaCra10}, Joshi and Shiu \cite{JoShiu13}, Banaji and Pantea \cite{BaPa16}, Conradi et al. \cite{CFMW17}. Under the restrictive assumption of mass action, see the works by Rendall and coauthors \cite{rendall1, rendall2, rendall3}, Dickenstein et al. \cite{Dickensteinetal19}, Shiu and de Wolff \cite{Shiu19}, Feliu et al. \cite{Timo20}.\\

One attractive mathematical possibility to detect multistationarity is to identify a saddle-node (SN) bifurcation. A bifurcation is a sudden qualitative change in the system behavior according to a small change in the parameter values. A saddle-node bifurcation occurs when two equilibria, e.g. one stable and one unstable, collide and disappear. Hence, the occurrence of such a bifurcation at a positive equilibrium automatically implies an area of parameters with at least two positive equilibria. Under the assumption of mass action, saddle-node bifurcations for biochemical systems have been {addressed} by Conradi et al. \cite{Conradi2007} and Domijan and Kirkilionis \cite{DomKirk09}. In both these contributions, the abstract conditions leading to the bifurcation have been reformulated in the polynomial language of mass action. Otero--Muras and coauthors used computational methods to detect saddle-node bifurcations in biochemical systems, see for example \cite{Otero18}.  Okada et al. \cite{Okada21} translated the bifurcation conditions from the Jacobian to an \emph{augmented} matrix, which allows them to confine a potential bifurcation behavior in a certain subnetwork. To the best of our knowledge, abstract network conditions that characterize saddle-node bifurcations have not yet been obtained for systems arising from chemical reaction networks. This is the focus of the present paper.\\

The standard saddle-node bifurcation theorem  for ordinary differential equations (ODEs) reads as follows.
\begin{theorem}[Saddle-node bifurcation,  \cite{GuHo84}]\label{SN}
Let $\dot{x}=g(x, \lambda)$ be an ODEs system in $\mathbb{R}^M$ depending on a single parameter $\lambda$. When $\lambda=\lambda^*$, assume that there is an equilibrium $\bar{x}$ for which the following hypotheses are satisfied:
\begin{enumerate}
\item[(SN1)] The Jacobian $G(\bar{x},\lambda^*):=D_x g(\bar{x},\lambda^*)$ has an algebraically simple eigenvalue 0 with right eigenvector $v$ and left eigenvector $w$.  $G(\bar{x},\lambda^*)$ has $\kappa$ eigenvalues with negative real part and $(M-\kappa-1)$ eigenvalues with positive real parts (counting multiplicity).
\item[(SN2)] $\langle w, \partial_\lambda g(\bar{x}, \lambda^*) \rangle  \neq 0$
\item[(SN3)] $w^T \; \partial_x^2 g(\bar{x}, \lambda^*)[v, v] \neq 0$
\end{enumerate}
Then there is a smooth curve of equilibria in $\mathbb{R}^M \times \mathbb{R}$ passing through $(\bar{x}, \lambda^*)$, tangent to the hyperplane $\mathbb{R}^M \times \{\lambda^*\}$. Depending on the signs of the expressions in (SN2) and (SN3), there are no equilibria near $(\bar{x}, \lambda^*)$ when $\lambda < \lambda^*$ ($\lambda > \lambda^*$) and two equilibria near $(\bar{x}, \lambda^*)$ for each parameter value $\lambda > \lambda^*$ ($\lambda <  \lambda^*$). The two equilibria for $\dot{x}=g(x, \lambda)$ near $(\bar{x}, \lambda^*)$ are hyperbolic and have stable manifolds of dimensions   $\kappa$ and $\kappa+1$, respectively. The set of equations $\dot{x}=g(x, \lambda)$ which satisfy (SN1)-(SN3) is open and dense in the space of $C^\infty$ one-parameter families of vector fields with an equilibrium at $(\bar{x},\lambda^*)$ with a zero eigenvalue.
\end{theorem}
Condition (SN1) is the necessary spectral condition: an algebraically simple eigenvalue zero of the Jacobian, at an equilibrium. Conditions (SN2) and (SN3) sufficiently guarantee the proper nonlinear unfolding of the bifurcation. We call \emph{degenerate saddle-node} the situation when conditions (SN1) and (SN2) are satisfied, but not (SN3). \emph{Nondegenerate saddle-node} indicates then the complete case where conditions (SN1)--(SN3) hold. Moreover, we will refer to properties that hold on an open and dense subset as \emph{generic}, albeit often in literature a generic set is more generally defined as a set of \emph{second Baire category} \cite{Baire}, i.e., a countable intersection of open and dense subsets. Theorem \ref{SN} is stated this way by Guckenheimer and Holmes \cite{GuHo84}, without an explicit proof, for which Vanderbauwhede \cite{vander89} is a reference. The genericity part has also been addressed and elaborated by Sotomayor \cite{soto73}.\\

A chemical reaction turns reactants into products. Several connected reactions constitute a chemical reaction network. We investigate which networks can sustain saddle-node bifurcations, and consequently multistationarity. More precisely, to any network $\mathbf{\Gamma}$ we associate the following ODEs dynamical system:
\begin{equation}\label{maineq}
\dot{x}= g(x) := S f(x),
\end{equation}
where $x(t)>0 \in \mathbb{R}^M$ is the vector of concentrations of the chemical species; the $M \times E$ matrix $S$ is the \emph{stoichiometric matrix}, the incidence matrix of the network; $f(x) \in \mathbb{R}^E$ is the vector of the reaction functions. We stress that we consider only strictly positive concentrations $x>0$: \emph{boundary equilibria}, where some of the concentrations $x_m$ are zero, fall beyond the scope of the present work. We address and answer the following question:
\begin{center}
\textit{For which networks $\mathbf{\Gamma}$ {does there exist} a choice of $f$\\ such that the associated dynamical system admits a saddle-node bifurcation?}
\end{center}

\vspace{0.2cm}

The precise form of $f$ is typically unknown in applications. Therefore, it is of great interest to obtain conditions only based on the network structure. Following this precise intention, we do not prescribe any specific form to $f$ but rather look into the entire set of functions satisfying only a few meaningful assumptions that make them reasonable as reaction functions, according to the following definition.
\begin{definition}[monotone chemical functions] \label{chemfeas} 
Let $j$ be a reaction and $f_j$ the associated reaction function. We call $f_j$ \emph{chemical} if 
\begin{enumerate}
\item $f_j$ depends only on the concentrations of the reactants of the reaction $j$;
\item $f_j$ is positive, i.e.,
$$f(x)>0,\text{ for every $x>0$}.$$
\end{enumerate}
We call a chemical function $f_j$ \emph{monotone} if
\begin{enumerate}
\setcounter{enumi}{2}
\item $f'_{jm}(x):=\frac{\partial f_j (x)}{\partial x_m} > 0$,  for any species $m$ reactant of $j$ and $x>0$.
\end{enumerate}
\end{definition}
Widely used and standard kinetic schemes as mass action \cite{HFJ72}, Michaelis--Menten \cite{MM13}, and Hill kinetics \cite{Hill10} follow Definition \ref{chemfeas}. {However, condition 1 excludes dependencies $f'_{jm}(x) \neq 0$ not expressed by the stoichiometry. Regulatory terms, i.e. $f'_{jm}(x) \neq 0$ with $m$ not a reactant to $j$, both in form of \emph{activators} $f'_{jm}(x) > 0$ and \emph{inhibitors} $f'_{jm}(x) < 0$, are not taken in account here. Condition 2 considers the reaction $j$ as \emph{irriversible}. As addressed in Section \ref{setting}, a \emph{reversible} process is treated in this setting as two opposite irriversible ones. Condition 3 excludes nonmonotone reaction rates as, for example, substrate inhibition.  Furthermore, condition 3 actually requires monotone \emph{increasing} functions, i.e. $f'_{jm}(x) > 0$, as this case is more relevant. Mathematically, we could develop analogous results with the monotone \emph{decreasing} condition: $f'_{jm}(x) < 0$. Yet, a small straightforward technicality must be taken in account: if the reaction functions are monotone increasing, then any product among nonzero partial derivatives is always positive. In contrast, for monotone decreasing reaction functions, the sign of the product depends on the number of factors, of course. We proceed assuming always monotone increasing functions with no further specification.} \\

We address the bifurcation conditions symbolically. For a related approach in bifurcation analysis on networks, see the work by Fiedler \cite{F19} that concerns global Hopf bifurcation. This symbolic strategy relates to the theory of \emph{jets} \cite{arnold2012}. Aiming at a self-contained presentation, we proceed from scratch. We call $\mathbf{r} \in \mathbb{R}^M_{>0}$ {the vector of equilibrium rates} that $f(\bar{x})$ attains at an equilibrium $\bar{x}$. The \emph{equilibrium constraints} define $\mathbf{r}$ and simply read
\begin{equation}\label{eqconst}
S\mathbf{r}=0, \quad\quad\quad \text{with }r_j>0\quad\text{for every $j$}.
\end{equation}
In particular, $\mathbf{r}$ is any positive right kernel vector of the stoichiometric matrix $S$. Throughout the paper, we only consider networks whose stoichiometric matrix $S$ admits a positive right kernel vector, i.e., admitting an equilibrium for a certain choice of chemical functions $f$. Without this basic assumption, addressing equilibria bifurcations would be meaningless. On the other hand, the \emph{bifurcation constraints} concern derivatives. We use the notation $$\mathbf{r}'=\{r'_{jm}\}_{j\in \mathbf{E}, \; m\in \mathbf{M}}$$ for the values, which the nonvanishing first derivatives {$f'_{jm}(\bar{x}):= \partial f_j (\bar{x}) / \partial x_m$} attain at the bifurcating equilibrium $\bar{x}$. Analogously, we use the notation $$\mathbf{r}''=\{r''_{jmn}\}_{j\in \mathbf{E}, \; m, n \in \mathbf{M}}$$ for the values of the second derivatives 
$$f''_{jmn}(\bar{x}):=\frac{\partial^2 f_j}{\partial{x_m}\partial{x_n}}(\bar{x}).$$ 

Firstly, we address symbolically conditions (SN1) and (SN2) in terms of the values $\mathbf{r}'$ alone. Secondly, we address the condition (SN3) in terms of the values $\mathbf{r}''$ alone. Finally, if conditions \eqref{eqconst} and (SN1)--(SN3) are satisfied by an independent choice ($\bar{\mathbf{r}}$, $\bar{\mathbf{r}}'$, $\bar{\mathbf{r}}''$), we find proper $f$ such that 
\begin{equation}\label{fconstr}
f(\bar{x})=\bar{\mathbf{r}},\quad\quad\quad f'_{jm}(\bar{x})=\bar{r}'_{jm} \quad\text{ and }\quad f''_{jmn}(\bar{x})=\bar{r}''_{jmn},
\end{equation}
for a positive equilibrium value $\bar{x}$ and any $j$, $m$, and $n$.\\

In this sense, we say that a chemical network $\mathbf{\Gamma}$ \emph{admits a saddle-node bifurcation} if there is a choice of $f(x, \lambda)$, within the class of monotone chemical functions, such that the assumptions of Theorem \ref{SN} hold. Of course, proving independently the conditions in terms of $(\mathbf{r}, \mathbf{r}', \mathbf{r}'')$ always implies the existence of a monotone chemical function $f$ for which all the bifurcation conditions (SN1)--(SN3) are satisfied, at any {choice of a} positive $\bar{x}$: the class of monotone chemical functions is clearly wide enough to include a nonlinearity $f$ satisfying \eqref{fconstr}.
However, even in its generosity, nature may not always provide us with such a freedom of choice, and typically given parametric class of functions (\emph{kinetics}) are used to model the reaction network. The validity of {the} results, when restricted to a certain kinetics must be further checked. In particular, we need the parametric freedom to assign independently the function value $\mathbf{r}$ and its first derivative value $\mathbf{r}'$, at least. We prove that this is possible in the parametric class of Michaelis-Menten kinetics. The only obstacle for multistationarity might reside in the tangency of the curve of equilibria, condition (SN3) of Theorem \ref{SN}, see example \ref{ex5}. The slightly more general Hill kinetics already provides the parametric freedom to conclude always a nondegenerate bifurcation result, in the present setting. On the contrary, polynomial mass action kinetics does not equally provide such parametric freedom. We show a mass-action example undergoing a saddle-node bifurcation in \ref{ex3}, to foster discussion.\\

We base the  results on the language of \emph{Child Selections}. A Child Selection $\mathbf{J}$ is an injective map associating to each species $m$ a reaction $j$, in which $m$ participates as a reactant, see Definition \ref{CSDEF}. The Jacobian determinant of the system, $\operatorname{det}G$, can be expanded along Child Selections (Proposition 2.1 of \cite{VGB20}) as:
\begin{equation}\label{detexp}
\operatorname{det}G=\sum_\mathbf{J} \alpha_{\mathbf{J}} \; \prod_{m \in \mathbf{M}} r'_{\mathbf{J}(m)m},
\end{equation}
where $\alpha_{\mathbf{J}}$ is a coefficient structurally associated to any Child Selection. Note that $\operatorname{det}G$ can then be interpreted as a multilinear homogenous polynomial $P(\mathbf{r'}):=\operatorname{det}G(\mathbf{r'})$, considering $\mathbf{r}'>0$ as independent real variables. Throughout, for simplicity of presentation, we assume the existence of at least one Child Selection $\mathbf{J}$ with $\alpha_\mathbf{J}\neq 0$, implying 
$$P(\mathbf{r'})\not\equiv0.$$
This excludes a permanent eigenvalue zero of $G$ and allows us to focus directly on solving $P(\mathbf{r'})=0$ without considering any reduced system. This assumption also excludes conserved linear combinations of the concentrations $x_m(t)$ for the whole network and may not be restrictive in itself: for instance, many metabolites in metabolic networks have a decay outflow reaction. The first main result, {discussed in Section \ref{main}}, characterizes the solvability of $P(\mathbf{r'})=0$ in terms of Child Selections.
\begin{thm*}
The multilinear homogeneous polynomial
$$P(\mathbf{r}'):=\operatorname{det} G(\mathbf{r}') $$
has a positive root $\bar{\mathbf{r}}'>0$ if and only if there exist two Child Selections $\mathbf{J}_1$, $\mathbf{J}_2$ such that $$\alpha_{\mathbf{J}_1}\alpha_{\mathbf{J}_2} < 0.$$
\end{thm*}
The above theorem characterizes the networks admitting a singular Jacobian. As stated in \ref{SN}, genericity of saddle-node bifurcations suggests that singular Jacobians indicate a nondegenerate bifurcation in most applications. However, \cite{V22} presents a ``pathological" network whose Jacobian $G$ possesses either no or multiple eigenvalue zero, for any choice of monotone chemical functions $f$. Even if rare and unexpected, such a case must be technically excluded. The second main result, Theorem \ref{SNmy}, provides a sufficient structural condition to have a saddle-node bifurcation. We define a saddle-node pair (SN-pair) of Child Selections, satisfying a further algebraic condition excluding multiple eigenvalues zero.  The presence of an SN-pair of Child Selections in the network guarantees the bifurcation behavior. Theorem \ref{SNmy} essentially reads: \emph{If the network possesses an SN-pair of Child Selections, then the network admits a saddle-node bifurcation.} The bifurcation parameter $\lambda$ is introduced parametrizing one single reaction function $f_j$, identified by an SN-pair of Child Selections.\\
 
The paper is organized as follows: Section \ref{setting} formalizes the mathematical setting, and Section \ref{CSPCS} introduces the language of Child Selections. The main results are presented in Section \ref{main}. Sections \ref{evzero}, \ref{almult}, \ref{unfolding}, and \ref{kinetics} build up the arguments needed to prove the main results. In particular, Section \ref{evzero} discusses networks possessing an eigenvalue zero; Section \ref{almult} addresses the multiplicity of such eigenvalue; Section \ref{unfolding} presents the unfolding of the bifurcation; Section \ref{kinetics} reads the results with explicit parameter choices for Hill and Michaelis-Menten kinetics, serving both as a specific example and as a general procedure on how to implement the results in given dynamical models. Section \ref{examples} lists four examples: \ref{ex1} a network motif giving rise to saddle-node bifurcation;  \ref{ex2} a saddle-node bifurcation identified in the central carbon metabolism of \emph{E.coli}; \ref{ex3} a mass-action example; \ref{ex5} an example of a network that admits only a degenerate saddle-node when endowed with Michaelis-Menten. Section \ref{discussion} concludes the paper with the discussion. Section \ref{proofs} lists all proofs.\\

\section{Setting}\label{setting}

A chemical reaction network $\mathbf{\Gamma}$ is a pair of sets $\{\mathbf{M},\mathbf{E}\}$: $\mathbf{M}$ is the set of chemical species or metabolites, and $\mathbf{E}$ is the set of reactions. Both sets are finite with cardinalities $|\mathbf{M}|=M$ and $|\mathbf{E}|=E$. Letters $m, n \in \mathbf{M}$ and $j, h \in \mathbf{E}$ refer to species and reactions, respectively. \\

A reaction $j$ is an ordered association of two positive linear combinations of species:
\begin{equation} \label{reactionj}
 j: \quad s^{j}_1m_1+...+s^{j}_Mm_M \underset{j}{\longrightarrow} \tilde{s}^{j}_1m_1+...+\tilde{s}^{j}_Mm_M.
\end{equation}
The nonnegative {real} coefficients $s^{j},\tilde{s}^j$ are called \emph{stoichiometric} coefficients. Chemical networks {often} deal with integer stoichiometric coefficients. The \emph{reactants} (resp., \emph{products}) of the reaction $j$ are the species appearing at the left (resp., right) of \eqref{reactionj} with nonzero stoichiometric coefficient. Chemical systems are often open systems: \emph{inflow reactions} are then reactions with no reactants ($s^j_m=0$ for every $m$) and \emph{outflow reactions} are reactions with no products ($\tilde{s}^j_m=0$ for every $m$). The $M \times E$ stoichiometric matrix $S$ is the matrix of all ordered stoichiometric coefficients:
\begin{equation}\label{smatrix}
S_{mj}:= \tilde{s}^j_{m} -s^j_m.
\end{equation}
This way we assign a fixed order to each reaction: we model a reversible reaction
\begin{equation*}
j: \quad A+2B \underset{j}{\rightleftharpoons} A+2C
\end{equation*}
simply as two irreversible reactions
\begin{equation} \label{j_1}
j_1: \quad A+2B \underset{j_1}{\longrightarrow} A+2C\quad \text{ and }\quad j_2: \quad A+2C \underset{j_2}{\longrightarrow} A+2B.
\end{equation}
We use the notation $S^j$ for the column of the stoichiometric matrix $S$ associated to the reaction $j$. For example, in a network of four species $\{A,B,C,D\}$, reaction $j_1$ in \eqref{j_1} is represented as the $j_1^{th}$ column of the stoichiometric matrix $S$ as
\begin{equation*}
S^{j_1}=
\begin{blockarray}{cc}
 & j_1 \\
\begin{block}{c(c)}
  A & 0\\
  B &  -2\\
  C & 2\\
  D & 0\\
\end{block}
\end{blockarray}\;.
\end{equation*}
Let $x \ge 0$ be the $M$-vector of chemical concentrations. Under the assumption that the reactor is well mixed, spatially homogeneous, and isothermal, the dynamics $x(t)$ of the concentrations satisfy the following system of ODEs:
\begin{equation} \tag{1}
\dot{x}=g(x):=S f(x),
\end{equation}
where $S$ is the $M \times E$ stoichiometric matrix \eqref{smatrix} and $f(x)$ is the $E$-vector of the \emph{reaction functions}. Without any reactant, we consider as constant the reaction function of inflow reactions $j_f$:
\begin{equation*}
f_{j_f}(x)\equiv F_{j_f}.
\end{equation*}
For any other reaction $j$, we only require that $f$ is \emph{monotone chemical}, as defined in \ref{chemfeas}.\\

\section{Child Selections and Partial Child Selections} \label{CSPCS}
We introduce the main tools. 
\begin{definition}[Child Selections \cite{BF18}] \label{CSDEF}
A \emph{Child Selection} is an injective map $\mathbf{J}: \textbf{M} \longrightarrow \textbf{E}$, which associates to every species $m \in \textbf{M}$ a reaction $j \in \textbf{E}$ such that $m$ is a reactant of reaction $j$.
\end{definition}
The notation $j \in \mathbf{J}$ indicates that there exists a species $m \in \mathbf{M}$ such that $j = \mathbf{J}(m)$. Let now $S^\mathbf{J}$ indicate the matrix whose $m^{th}$ column is the $\mathbf{J}(m)^{th}$ column of $S$. In particular, the columns of $S^\mathbf{J}$ correspond one-to-one and following the order to the reactions
$$\mathbf{J}(m_1), \; \mathbf{J}(m_2), \;...\; , \;\mathbf{J}(m_{M-1}), \;\mathbf{J}(m_{M}).$$
We associate to each Child Selection $\mathbf{J}$ the coefficient
$$\alpha_\mathbf{J} : = \operatorname{det} S^\mathbf{J}.$$
Let $G:=\partial_x g(x)$ indicate the Jacobian matrix of \eqref{maineq}. The determinant of $G$ can be expressed in terms of Child Selections \cite{BF18,VGB20}:
\begin{equation} \label{vita}
\operatorname{det}G=\sum_\mathbf{J}  \alpha_\mathbf{J}  \cdot \prod_{m\in \mathbf{M}} f'_{\mathbf{J}(m)m}(x),
\end{equation}
The sum runs on all Child Selections. Let us consider the matrix $G$ symbolically, i.e., we consider the nonvanishing partial derivatives $f'_{\mathbf{J}(m)m}$ as independent positive variables $r'_{jm}:=f'_{\mathbf{J}(m)m}$. We then interpret \eqref{vita} as a multilinear homogenous polynomial of order $M$, in the variables $\mathbf{r}'$. The notation $\mathbf{r}'[\mathbf{J}]$ indicates the monomial of the variables $r'_{jm}=r'_{\mathbf{J}(m)m}$. That is,
$$\mathbf{r}'[\mathbf{J}]:=\prod_{m\in \mathbf{M}} r'_{\mathbf{J}(m)m}.$$
In this sense,
\begin{equation}\label{vita2}
P(\mathbf{r}'):=\operatorname{det}G(\mathbf{r}')=\sum_\mathbf{J}  \alpha_\mathbf{J}  \mathbf{r}[\mathbf{J}]
\end{equation}
is the symbolic version of polynomial \eqref{vita}.\\

We call $\alpha_\mathbf{J}$ the \emph{behavior coefficient}. Depending on the sign of $\alpha_{\mathbf{J}}$ we classify a Child Selection as follows. We call a Child Selection $\mathbf{J}$ \emph{zero} if $\alpha_\mathbf{J}= 0$. On the contrary, we call $\mathbf{J}$ a \emph{nonzero} Child Selection if $\alpha_\mathbf{J}\neq0$. In the latter case, we say that $\mathbf{J}$ is \emph{good} if $\operatorname{sign}\alpha_\mathbf{J}=(-1)^M$, and \emph{bad} otherwise. To clarify the naming, let us consider a system that possesses {a single stable equilibrium} for any choice of $f$. This requires the Jacobian of such equilibrium to have either only eigenvalues with negative real part or pairs of purely imaginary complex conjugated eigenvalues, and it excludes saddle-node bifurcations, of course. Assuming at least {one} eigenvalue with negative real part, the sign of a nonsingular Jacobian is
$$\operatorname{sign}\operatorname{det}SR= (-1)^M.$$
Via \eqref{vita2}, this ``stable'' sign is automatically implied if there are no bad Child Selections. In the opposite direction, a loss of stability of an equilibrium {via a sign-change of a single eigenvalue} necessarily implies the existence of at least one bad Child Selection. Furthermore, any Child Selection naturally identifies a subnetwork consisting only of reactions $j \in \mathbf{J}(\mathbf{M})$. In \cite{VGB20}, the behavior of any Child Selection has been structurally characterized. That analysis showed that certain classes of Child Selections, commonly found in metabolic networks, are always good. As a consequence, we observe a clear predominance of good Child Selections in metabolic networks, whereas the few bad Child Selections hint at stability change and bifurcations.\\

A natural distance can be assigned to the set of Child Selections.
\begin{definition}[Distance of Child Selections]
Let $\mathbf{J}_1, \mathbf{J}_2$ be two Child Selections. The distance $d(\mathbf{J}_1, \mathbf{J}_2)$ is the number of species $m \in \mathbf{M}$ such that $\mathbf{J}_1(m) \neq \mathbf{J}_2(m)$.
\end{definition}
\begin{remark}
In literature, this distance is often called \emph{Hamming} distance \cite{nowak06}. 
\end{remark}
This distance admits a natural concept of minimality, in the following sense.
\begin{definition}[Minimal distance]
Let $\mathbf{J}_1$, $\mathbf{J}_2$ be two nonzero Child Selections at distance $d=\delta$. We say that $\mathbf{J}_1$ and $\mathbf{J}_2$ are at \emph{minimal distance} if any Child Selection $\mathbf{J}_3$ such that 
\begin{equation}\label{mindisteq}
\begin{cases}
d(\mathbf{J}_1, \mathbf{J}_3) < \delta\\
d(\mathbf{J}_2, \mathbf{J}_3) < \delta\\
\end{cases}
\end{equation}
is zero, i.e., $\alpha_{\mathbf{J}_3} = 0.$
\end{definition}
\begin{remark}
The above definition always applies if $\mathbf{J}_1$ and $\mathbf{J}_2$ are at distance $d=1$, since no Child Selection $\mathbf{J}_3$ satisfies \eqref{mindisteq}.
\end{remark}

We conclude this section with a related concept: the \textit{Partial Child Selections}.
\begin{definition} [Partial Child Selections]
A \emph{Partial Child Selection} $\mathbf{J}^{\vee m}$ is an injective map:
$$\mathbf{J}^{\vee m} : \mathbf{M}\setminus \{m\} \longrightarrow \mathbf{E},$$
associating to each species $n \neq m$ a reaction $j$ such that $n$ is a reactant of $j$.\\
\end{definition}
Let us pick a metabolite $m_i$ and consider $1,..., i, ... M$ without loss of generality. In analogy to the submatrix $S^{\mathbf{J}}$ for a Child Selection $\mathbf{J}$, the expression $S^{\mathbf{J}^{\vee m_i}}$ indicates the $M \times (M-1)$ matrix with columns corresponding one-to-one, and following the order, to the reactions
$$\mathbf{J}^{\vee m_i}(m_1),\;...\;,\;\mathbf{J}^{\vee m_i}(m_{i-1}), \;\mathbf{J}^{\vee m_i}(m_{i+1}),\;...\;,\; \mathbf{J}^{\vee m_i}(m_M).$$
The first column is the stoichiometric column $S^{j_1}$ of the reaction $j_1=\mathbf{J}^{\vee m_i}(m_1)$ and the $i^{th}$ column is the stoichiometric column $S^{j_i}$ of the reaction $j_i=\mathbf{J}^{\vee m_i}(m_{i+1})$, and so on. We associate to each Partial Child Selection $\mathbf{J}^{\vee m}$ the behavior coefficient
$$\beta_{\mathbf{J}^{\vee m}} : =  \operatorname{det}  S^{\mathbf{J}^{\vee m}}_{\vee m},$$
where the notation $S^{\mathbf{J}^{\vee m}}_{\vee m}$ indicates the $(M-1)\times(M-1)$ matrix obtained from $S^{\mathbf{J}^{\vee m}}$ by removing the $m^{th}$ row. If the behavior coefficient $\beta_{\mathbf{J}^{\vee m}}$ is zero (resp., nonzero) we call the Partial Child Selection $\mathbf{J}^{\vee m}$ \emph{zero} (resp., \emph{nonzero}), accordingly.

\section{Main results}\label{main}

The first result is a characterization of networks that admit a singular Jacobian.

\begin{theorem}\label{det0}
The multilinear homogeneous polynomial
$$P(\mathbf{r}'):=\operatorname{det} G(\mathbf{r}') $$
has a positive root $\bar{\mathbf{r}}'>0$ if and only if there exist two Child Selections $\mathbf{J}_1$, $\mathbf{J}_2$ such that $$\alpha_{\mathbf{J}_1}\alpha_{\mathbf{J}_2} < 0.$$
\end{theorem}
Theorem \ref{det0} is the structural characterization of a necessary spectral condition for a saddle-node bifurcation. However, the existence of two Child Selections $\mathbf{J}_1$, $\mathbf{J}_2$ with $\alpha_{\mathbf{J}_1}\alpha_{\mathbf{J}_2} < 0$ does not guarantee that there exists a positive root $\bar{\mathbf{r}}'$ of $P(\mathbf{r}')$ such that the associated Jacobian $G(\bar{\mathbf{r}}')$ has an \emph{algebraically simple} eigenvalue zero. See \cite{V22} for a counterexample. For the simplicity, we need a further condition. 
\begin{definition}[SN-pair of Child Selections]\label{snpairs}
We call two nonzero Child Selections $\mathbf{J}_1$, $\mathbf{J}_2$ a \emph{saddle-node pair} (SN-pair) if the following conditions all hold true:
\begin{enumerate}
\item $\mathbf{J}_1$ and $\mathbf{J}_2$ are at minimal distance;
\item $\alpha_{\mathbf{J}_1} \alpha_{\mathbf{J}_2} < 0$;
\item there exists a species $\tilde{m}$ with $\mathbf{J}_1(\tilde{m}) \neq \mathbf{J}_2(\tilde{m})$ and a nonzero Partial Child Selection $\mathbf{J}^{\vee \tilde{m}}$ such that $\mathbf{J}^{\vee \tilde{m}}(n)=\mathbf{J}_1(n)$ or $\mathbf{J}^{\vee \tilde{m}}(n)=\mathbf{J}_2(n)$ for every $n \neq \tilde{m}$.
\end{enumerate}
\end{definition}

We can now state the main result of this paper.

\begin{theorem}\label{SNmy}
Assume that the network possesses an SN-pair of Child Selections $\mathbf{J}_1, \mathbf{J}_2$. Then there exists a choice of monotone chemical functions $f$ such that the associated dynamical system
$$\dot{x}=S f(x, \lambda)$$
undergoes a saddle-node bifurcation at a positive equilibrium $\bar{x}$ for a bifurcation value $\lambda^*$. The bifurcation parameter $\lambda$ parametrizes the function $f_{\eta}$ of a reaction $\eta$ such that $\mathbf{J}_1(m^*)=\eta \neq \mathbf{J}_2(m^*)$, for a species $m^*$. 
\end{theorem}

Theorem \ref{SNmy} states that the existence of an SN-pair of Child Selections is a sufficient condition for the network to admit a saddle-node bifurcation. We derive the following corollary that provides a method to identify bifurcation motifs in reaction networks.

\begin{corollary}\label{mindist1}
Let a network $\mathbf{\Gamma}=(\mathbf{M},\mathbf{E})$ possess an SN-pair of Child Selections $(\mathbf{J}_1, \mathbf{J}_2)$ and let $\tilde{\mathbf{\Gamma}}=(\tilde{\mathbf{M}}, \tilde{\mathbf{E}})$ be a network having $\mathbf{\Gamma}$ as a subnetwork: i.e., $\mathbf{M}\subseteq \tilde{\mathbf{M}}$, $\mathbf{E} \subseteq \tilde{\mathbf{E}}$. Assume there exists a pair of nonzero Child Selections of $\tilde{\mathbf{\Gamma}}$, $(\tilde{\mathbf{J}}_1, \tilde{\mathbf{J}}_2)$, at minimal distance in $\tilde{\mathbf{\Gamma}}$, such that
\begin{equation*}
\begin{cases}
\tilde{\mathbf{J}}_1(m)=\mathbf{J}_1(m)\\
\tilde{\mathbf{J}}_2(m)=\mathbf{J}_2(m)\\
\end{cases},
\end{equation*}
for every $m \in \mathbf{M}$, and $\tilde{\mathbf{J}}_1(m)=\tilde{\mathbf{J}}_2(m)$ for every $m \in \tilde{\mathbf{M}} \setminus \mathbf{M}$. Then $(\tilde{\mathbf{J}}_1, \tilde{\mathbf{J}}_2)$ form an SN-pair of Child Selections of $\tilde{\mathbf{\Gamma}}$ and, in particular, $\tilde{\mathbf{\Gamma}}$ admits a saddle-node bifurcation.
\end{corollary}

\begin{remark}
In the case of minimal distance $d=1$, the minimality in the subnetwork $\mathbf{\Gamma}$ is always inherited by $\tilde{\mathbf{\Gamma}}$. 
\end{remark}

In the following sections, we will present all the arguments {that} prove the results.

\section{Eigenvalues zero}\label{evzero}

We solve $P(\mathbf{r}')=0$ by considering a rescaling of the variables $\mathbf{r}'$ with the introduction of a further parameter $\varepsilon >0$. The rescaling identifies two ``leading" monomials corresponding to two Child Selections at minimal distance.

\begin{lemma}\label{mindistsol}
Let $\mathbf{J}_1$ and $\mathbf{J}_2$ be two Child Selections at minimal distance. Then there exists an $\varepsilon$-rescaling of $\mathbf{r}'$ such that
\begin{equation}\label{lemmind}
P(\mathbf{r}')=\alpha_{\mathbf{J}_1} \mathbf{r}'[\mathbf{J}_1]+\alpha_{\mathbf{J}_2} \mathbf{r}'[\mathbf{J}_2]+ q(\varepsilon),
\end{equation}
where $q(\varepsilon)$ is a polynomial with $q(0)=0$.
\end{lemma}

Throughout the paper, we extensively refer to such $\varepsilon$-rescaling. If $\alpha_{\mathbf{J}_1}\alpha_{\mathbf{J}_2}<0$,  $$\alpha_{\mathbf{J}_1} \mathbf{r}'[\mathbf{J}_1]=-\alpha_{\mathbf{J}_2} \mathbf{r}'[\mathbf{J}_2]$$ is a positive solution of $P=0$ at the limit $\varepsilon=0$. In the proof of Theorem \ref{det0}, we employ the implicit function theorem to extend this solution to positive $\varepsilon>0$. The bridge between Lemma \ref{mindistsol} and Theorem \ref{det0} is then provided by the following lemma.

\begin{lemma} \label{mainlemma}
There exist two Child Selections $\mathbf{J}_1$ and $\mathbf{J}_2$ with 
$\alpha_\mathbf{J_1} \alpha_\mathbf{J_2} < 0$ 
if and only if there exist two Child Selections $\mathbf{J}_3$ and $\mathbf{J}_4$ \emph{at minimal distance} with 
$\alpha_\mathbf{J_3} \alpha_\mathbf{J_4} < 0.$
\end{lemma}

Lemma \ref{mainlemma} is inspired by Balinski's theorem on convex polyhedra \cite{Balinski61}. Informally: consider the Newton polytope $N$ generated by the determinant polynomial $P(\mathbf{r}')$. By convexity, an expansion of the type \eqref{lemmind} is possible if and only if there is an edge $e$ of $N$ adjacent to the vertices associated to $\mathbf{J}_1$ and $\mathbf{J}_2$. Moreover, the expansion \eqref{lemmind} provides a positive root of $P(\mathbf{r}')$ if and only if the coefficient sign of adjacent vertices of $e$ is opposite. There is such an edge if and only if there are vertices of different sign: in fact, Balinski's theorem states the connectedness of the graph associated to $N$; hence it is enough to consider any path joining two vertices of different sign and find the first sign-switch. Even though this context is inspirational for the results, we proceed in a more elementary setting and we do not require any knowledge of convex polytopes.\\

\section{Algebraic multiplicity} \label{almult}

Theorem \ref{det0} characterizes a network for which the Jacobian $G$ of the associated system \eqref{maineq} admits a zero eigenvalue. To address the multiplicity of such an eigenvalue zero, we study the \emph{adjugate matrix} (transpose of the cofactor matrix) of $G$, $\operatorname{Adj}G$. We recall two propositions from \cite{V22}.

\begin{proposition}[\cite{V22}]
The Jacobian $G(\mathbf{r}')$ has an algebraically simple eigenvalue zero at $\bar{\mathbf{r}}'>0$ if and only if
\begin{equation} \label{snsyst}
\begin{cases}
P(\bar{\mathbf{r}}'):=\operatorname{det}G(\bar{\mathbf{r}}')=0;\\
A(\bar{\mathbf{r}}'):=\operatorname{tr}\operatorname{Adj}G(\bar{\mathbf{r}}')\neq 0.\\
\end{cases}.
\end{equation}
\end{proposition}

In analogy to the expansion \eqref{vita2} for $P(\mathbf{r}')$, the polynomial $A(\mathbf{r}')$ can be expanded along Partial Child Selections.

\begin{proposition}[\cite{V22}] \label{AdG}
Let $G$ be the Jacobian matrix of the system \eqref{maineq} and let $\operatorname{Adj}G^m_m$ indicate the $m^{th}$ diagonal entry of its adjugate. Then the following expansion holds:
\begin{equation} \label{23}
\operatorname{Adj}G^m_m(\mathbf{r}')= \sum_{\mathbf{J}^{\vee m}} \; \beta_{\mathbf{J}^{\vee m}} \;\mathbf{r}'[\mathbf{J}^{\vee m}],
\end{equation}
where $\mathbf{J}^{\vee m}$ are Partial Child Selections and the notation $\mathbf{r}'[\mathbf{J}^{\vee m}]$ indicates the multilinear monomial of degree $M-1$:
$$\mathbf{r}'[\mathbf{J}^{\vee m}]=\prod_{n \neq m} r'_{\mathbf{J}^{\vee m}(n)n}.$$
In particular,
\begin{equation*}
A(\mathbf{r}'):=\operatorname{tr}\operatorname{Adj}G(\mathbf{r}')=\sum_{m\in \mathbf{M}} \sum_{\mathbf{J}^{\vee m}} \; \beta_{\mathbf{J}^{\vee m}} \;\mathbf{r}'[\mathbf{J}^{\vee m}].
\end{equation*}
\end{proposition}

The paper \cite{V22} also presents a degenerate network for which $$\operatorname{det}G(\bar{\mathbf{r}}')=0\quad\quad\quad \Rightarrow \quad\quad\quad \operatorname{tr}\operatorname{Adj}G(\bar{\mathbf{r}}')= 0,$$
and hence the system \eqref{snsyst} is never satisfied. The presence of an SN-pair of Child Selection excludes this degeneracy, as the following lemma states.
\begin{lemma}\label{simple0}
Assume that the network possesses an SN-pair of Child Selections. Then the polynomial system \eqref{snsyst} has a positive solution $\bar{\mathbf{r}}'>0$. In particular, at $\bar{\mathbf{r}}'$ the Jacobian $G(\bar{\mathbf{r}}')$ possesses an algebraically simple eigenvalue zero.
\end{lemma}

\section{Nonlinear unfolding}\label{unfolding}

Let us consider a network with an SN-pair of Child Selections $\mathbf{J}_1$ and $\mathbf{J}_2$. Lemma \ref{simple0} guarantees the existence of a positive choice $\bar{\mathbf{r}}'>0$, such that the Jacobian $G(\bar{\mathbf{r}}')$ possesses an algebraically simple eigenvalue zero, i.e., spectral condition (SN1) of \ref{SN}. The nonlinear unfolding comprises conditions (SN2) and (SN3). Let $\eta$ be a reaction such that $\mathbf{J}_1(m^*)= \eta \neq j_2= \mathbf{J}_2(m^*)$, for a species $m^*$. We unfold the bifurcation by a {$\lambda$-parametrization} of the reaction $\eta$. In particular, the bifurcation parameter $\lambda$ appears in the reaction function $f_\eta$, only. To the bifurcation point ($\bar{x}, \lambda^*$) corresponds the bifurcation value $\bar{\mathbf{r}}'$. The nondegeneracy condition (SN2) requires that the derivative of the vector field with respect to $\lambda$ is not in the range of the Jacobian at the bifurcation point:
$$\langle w,\partial_\lambda g(\bar{x}, \lambda^*)\rangle  \neq 0,$$
for $w$ left eigenvector of the Jacobian $G(\bar{\mathbf{r}}')$. The first consequence of our parametrization choice is that the vector $\partial_\lambda g$ is parallel to the stoichiometric vector $S^{\eta}$ of reaction $\eta$ and 
$$\langle w,\partial_\lambda g(\bar{x}, \lambda^*)\rangle  \neq 0\quad \quad \quad \Leftrightarrow \quad \quad \quad \langle w,S^{\eta } \rangle \neq 0.$$
We have the following lemma.

\begin{lemma}[SN2]\label{SN2}
Assume that the network possesses an SN-pair of Child Selections $\mathbf{J}_1$ and $\mathbf{J}_2$. Let $\eta$ be a reaction such that $\mathbf{J}_1(m^*)= \eta \neq j_2= \mathbf{J}_2(m^*)$. Then there exists a positive root $\bar{\mathbf{r}}'>0$ of \eqref{snsyst} such that 
$$\langle w, S^{\eta } \rangle \neq 0,$$
where $w$ is a left kernel vector of the Jacobian $G(\bar{\mathbf{r}}')$ and $S^{\eta }$ is the stoichiometric column of reaction $\eta$.
\end{lemma}

Lemma \ref{SN2} shows that condition (SN2) is always satisfied by our choice of $f(x, \lambda)$ at the bifurcation point $(\bar{x}, \lambda^*)$. The last step is discussing the tangency of the curve of equilibria at the bifurcation point. Condition (SN3) states that a quadratic tangency is sufficient. 

\begin{lemma}[SN3]\label{sn3noneq0}
Let $\bar{\mathbf{r}}'$ be any positive root of the system \eqref{snsyst}, with $w$ and $v$ respectively left and right kernel vectors of the Jacobian $G(\bar{\mathbf{r}}')$. Then
$$w^T \; \partial^2_x g(\bar{x},\lambda^*) [v,v] \not\equiv 0,$$ 
as a function of the second derivatives $\mathbf{r}''=f''(\bar{x})$. 
\end{lemma}

Lemma \ref{sn3noneq0} concludes that the system admits a saddle-node bifurcation, if we have enough parametric freedom to assign $\mathbf{r}''$ freely and independently from $\mathbf{r}$ and $\mathbf{r}'$. We discuss this in detail in Section \ref{kinetics}, where we present the applicability of the results for two kinetics of interest.

\section{Michaelis-Menten and Hill kinetics} \label{kinetics}

In this section, we apply the results to given kinetics. We discuss Hill kinetics, as a general mathematical form that comprises also Michaelis-Menten and mass action kinetics as particular cases. Hill kinetics is a relevant example of monotone chemical functions. The mathematical form of a reaction $j$ according to Hill is:
\begin{equation}\label{MM}
f_j(x):=a_j\prod_{m\in\mathbf{M}} \Bigg( \frac{x_m^{c^j_m}}{(1+b^j_m x_m^{c^j_m})}\Bigg)^{s^j_m},
\end{equation}
where $s^j_m$ is the stoichiometric coefficient of species $m$ as reactant of the reaction $j$, and $a_j, b^j_m, c^j_m$ are positive parameters. Typically, $a_j, b^j_m$ are real, while $c^j_m$ is an integer, though irrelevant for the present mathematical description. Michaelis-Menten kinetics fixes $c^j_m = 1$ for all $j, m$. Mass action kinetics is recovered by considering the limit case $b^j_m=0, c^j_m = 1$ for all $j,m$. We write $\mathbf{a}$ to refer to the set of parameters $a_j$ for all reactions $j$. Analogously, we write $\mathbf{b}$ (resp., $\mathbf{c}$) for the set of parameters $b^j_m$ (resp., $c^j_m$), for all $j$ and $m$.\\

The results of this section can be summarized as follows: at any concentration value $\bar{x}$, the parametric freedom of Michaelis-Menten allows us to consider the values $\mathbf{r}$ of the function $f$, and the values $\mathbf{r}'$ of their derivatives $f'$ as independent parameters, via a careful choice of parameters $\mathbf{a}, \mathbf{b}$. Contrarily, the value of the second derivatives $\mathbf{r}''$ cannot be independently chosen. As a consequence, under the assumptions of Theorem \ref{SNmy}, we can always conclude that the network endowed with Michaelis-Menten kinetics possesses a positive equilibrium satisfying conditions (SN1) and (SN2) of Theorem \ref{SN}.  This is presented in Theorem \ref{S0MM}. However, we may never be able to find parameters that jointly satisfy also (SN3) of \ref{SN}. We present in Example \ref{ex5} a network showing such degeneracy. Theorem \ref{SNMM} provides then a sufficient condition to exclude this degeneracy in a Michaelis-Menten system. The degeneracy can always be avoided in the more general class of Hill kinetics, by a proper choice of the further available parameters $\mathbf{c}$, Theorem \ref{SNH}. {In contrast,} our network assumptions do not conclude a bifurcation result in the case of mass action kinetics. Nevertheless, Example \ref{ex3} presents a mass-action system undergoing a saddle-node bifurcation, where the construction is inspired by {the present} results.\\

We first present the two theorems for Michaelis-Menten kinetics, hence fixing $c^j_m=1$ for all reactions $j$ and species $m$ in the nonlinearity \eqref{MM}. We are thus left only with the choice of $\mathbf{a}, \mathbf{b}$.

\begin{theorem} \label{S0MM}
Assume that the network possesses an SN-pair of Child Selections $\mathbf{J}_1, \mathbf{J}_2$. Let $\eta$ be a reaction such that $\mathbf{J}_1(m^*)=\eta \neq \mathbf{J}_2(m^*)$, for the species $m^*$. Choose as bifurcation parameter $\lambda:=b^\eta_{m^*}$. Then there exists a choice of $\mathbf{a},\mathbf{b}$ such that the Michaelis-Menten system admits a positive equilibrium satisfying conditions (SN1) and (SN2) of Theorem \ref{SN}. 
\end {theorem} 
Michaelis-Menten kinetics does not guarantee a parameter choice such that the curve of equilibria at the bifurcation point has a quadratic tangency (SN3). For this reason, we present a condition that characterizes the nondegeneracy of a saddle-node bifurcation under Michaelis-Menten kinetics for the case where the SN-pair of Child Selections is at distance $d=1$. 
\begin{theorem} \label{SNMM}
Let $(\mathbf{J}_1, \mathbf{J}_2)$ be an SN-pair of Child Selections at distance $d=1$. Let $m^*$ be the unique species such that $\mathbf{J}_1(m^*)= \eta \neq j_2 =\mathbf{J}_2(m^*)$. Choose as bifurcation parameter  $\lambda=b^\eta_{m^*}$. Assume the following condition holds:
\begin{equation}\label{MMcond}
\frac{\alpha_{\mathbf{J}_2}}{\bar{r}_\eta}\bigg(1+\frac{1}{s^{\eta}_{m^*}}\bigg) \neq - \frac{\alpha_{\mathbf{J}_1}}{\bar{r}_{j_2}}\bigg(1+\frac{1}{s^{j_2}_{m^*}}\bigg),
\end{equation}
where $\bar{r}_\eta$ and $\bar{r}_{j_2}$ indicate the equilibrium rates \eqref{eqconst} relative to reaction $\eta$ and $j_2$, respectively. Then, there exists a choice of $\mathbf{a},\mathbf{b}$ such that the Michaelis-Menten system undergoes a saddle-node bifurcation according to the parameter $\lambda$.
\end{theorem}

In particular, Theorem \ref{SNMM} states that the degeneracy of the saddle-node depends on the ratio $\bar{r}_\eta/\bar{r}_{j_2}$, which is not uniquely fixed in most applications.  The degeneracy can be thus most often avoided by a proper choice of the equilibrium rates $\mathbf{r}$. See again Example \ref{ex5}. {Child Selections at greater minimal distance $d>1$ possess a quite special structure that will be addressed and described in a future publication. Such structure indicates also that the case $d=1$ is the most relevant, as it is the most likely to occur. See also the related discussion in Section \ref{discussion}. A general version of Theorem \ref{SNMM} for minimal distance $d>1$ requires the understanding of such structure, which exceeds the purposes of the present paper, and it is thus not addressed here.}\\

For the more general Hill kinetics, we can choose also parameters $\mathbf{c}\neq \mathbf{1}$. The result reads as follows.
\begin{theorem}\label{SNH}
Assume that the network possesses an SN-pair of Child Selections $\mathbf{J}_1, \mathbf{J}_2$. Let $\eta$ be a reaction such that $\mathbf{J}_1(m^*)=\eta \neq \mathbf{J}_2(m^*)$, for the species $m^*$. Choose as bifurcation parameter $\lambda:=b^\eta_{m^*}$. Then there exists a choice of $\mathbf{a}, \mathbf{b}, \mathbf{c}$ such that the Hill system undergoes a saddle-node bifurcation according to the parameter $\lambda$.
\end{theorem} 

Let us be explicit in the parameter choice: assume there exist positive $\bar{x}$, $\bar{\mathbf{r}}$, $\bar{\mathbf{r}}'$ such that:
\begin{equation}\label{posb1}
\begin{cases}
S \bar{\mathbf{r}} \mathbf = 0;\\
\operatorname{det}G(\bar{\mathbf{r}}')=0;\\
\frac{\bar{r}_j}{\bar{r}'_{jm}} \ge \frac{\bar{x}_m}{s^j_m} \quad \quad \quad \text{for every reaction $j$ and species $m$}.
\end{cases}
\end{equation}
Note that the three constraints \eqref{posb1} can be always satisfied for a network admitting a choice $\bar{\mathbf{r}}$, $\bar{\mathbf{r}}'$ satisfying the first two constraints: the third constraint follows by choosing big enough equilibrium flux $\bar{\mathbf{r}}$. We fix
\begin{equation}\label{bconst}
0<\mathpzc{b}_m^j:=\bigg(\frac{\bar{r}_j}{\bar{r}'_{jm}}\frac{s^j_m}{\bar{x}_m}-1\bigg)\frac{1}{\bar{x}_m},
\end{equation}
and
\begin{equation}\label{a}
\mathpzc{a}_j := \bar{r}_j \prod_{m\in\mathbf{M}} \Bigg( \frac{\bar{x}_m^{c^j_m}}{(1+\mathpzc{b}^j_m \bar{x}_m^{c^j_m})}\Bigg)^{-s^j_m}.
\end{equation}
A straightforward computation shows that the Hill function
$$f_j(x_m):=\mathpzc{a}_j  \prod_{ m}\Bigg( \frac{x_m^{c^j_m}}{(1+\mathpzc{b}^j_m x_m^{c^j_m})}\Bigg)^{s^j_m}$$
satisfies
\begin{equation*}
\begin{cases}
f_j(\bar{x}_m)=\bar{r}_j;\\
f_{jm}(\bar{x}_m)=\bar{r}'_{jm}.\\
\end{cases}.
\end{equation*}
Note that \eqref{bconst} and \eqref{a} do not require a fixed choice of $\mathbf{c}$ and hence hold true also for Michaelis-Menten, i.e. $\mathbf{c}=\mathbf{1}$. Furthermore, parameters $\mathbf{c}$ can be used to nudge 
$$ w^T \; \partial^2_x g(\bar{x},\lambda^*) [v,v]$$ away from the degenerate value 0, in the Hill case. We discuss it in detail in the proof of Theorem \ref{SNH}.

\section{Examples}\label{examples}
\subsection{Example I: Reversible feedback cycles} \label{ex1}
We present a family of networks admitting saddle-node bifurcations. Consider a reversible feedback cycle of length $M$:
$$m_1 \quad \overset{1}{\underset{4}{\myrightleftarrows{\rule{0.5cm}{0cm}}}} \quad m_2 \quad \overset{3}{\underset{6}{\myrightleftarrows{\rule{0.5cm}{0cm}}}} \quad ...\quad \overset{2M-3}{\underset{2M}{\myrightleftarrows{\rule{0.5cm}{0cm}}}} \quad m_M \quad \overset{2M-1}{\underset{2}{\myrightleftarrows{\rule{0.5cm}{0cm}}}} \quad 2m_1.$$
The feedback cycles generalize autocatalytic processes: walking along the cycle from left to right, one single molecule of $m_1$ produces two molecules of $m_1$, while from right to left, two molecules of $m_1$ reduce to one single molecule of $m_1$. We show that such a structure admits saddle-node bifurcations. More specifically, we can identify 2M different parameters triggering a saddle-node bifurcation. The system of $M$ differential equations reads:
\begin{equation*}
\begin{cases}
\dot{x}_1=-r_1(x_1)-2r_2(x_1)+r_4(x_2) + 2r_{2M-1}(x_M);\\
\dot{x}_i=-r_{2i-1}(x_i)-r_{2i}(x_i)+r_{2i+2}(x_{i+1})+r_{2i-3}(x_{i-1}),\quad \quad \text{for $i=2,...,M-1$};\\
\dot{x}_M=-r_{2M-1}(x_m)-r_{2M}(x_M)+2r_2(x_1)+r_{2M-3}(x_{M-1}).\\
\end{cases}
\end{equation*}
where $x_i$ is the concentration of $m_i$. An equilibrium is given by 
$$r_j \equiv \bar{r} \in \mathbb{R}_{>0}, \quad \quad \quad \text{for every $j$}.$$
There are only two nonzero Child Selections:
\begin{align*}
\begin{cases}
\mathbf{J}_1(x_i)=r_{2i-1}, \quad \quad \quad &\text{for every $i=1,...,M$};\\
\mathbf{J}_2(x_i)=r_{2i},\quad \quad \quad \;\;\; &\text{for every $i=1,...,M$}.
\end{cases}
\end{align*}
Since there are no other nonzero Child Selections, $\mathbf{J}_1$ and $\mathbf{J}_2$ are obviously at minimal distance $d(\mathbf{J}_1,\mathbf{J}_2)=M$. The behavior coefficients are opposite: $\alpha_{\mathbf{J}_1}=(-1)^{M-1}$ and $\alpha_{\mathbf{J}_2}=(-1)^{M}$, thus 
$$\operatorname{det}G(\mathbf{r}')=(-1)^{M-1}\mathbf{r}'[\mathbf{J}_1]+(-1)^{M}\mathbf{r}'[\mathbf{J}_2].$$
The determinant is zero if and only if $\mathbf{r}'[\mathbf{J}_1]=\mathbf{r}'[\mathbf{J}_2]$. Any Partial Child Selection
$$\mathbf{J}^{\vee m}(\mathbf{M}\setminus \{ m \})=\mathbf{J}_i(\mathbf{M}\setminus \{ m \}),$$
with $i=1$ or $i=2$, is nonzero. Hence $\mathbf{J}_1$ and $\mathbf{J}_2$ form an SN-pair of Child Selections, and the system admits a saddle-node bifurcation, via Theorem \ref{SNmy}.\\

To exemplify further, we compute all conditions explicitly under the assumption of Michaelis-Menten kinetics. We operate as described in Section \ref{kinetics} choosing arbitrary values. We fix $\bar{x}_i=1$ for every $i$, the values $\bar{r}_j=4$ for every $j$, and $\bar{r}_{jm}=1$ for every $j$ and $m$. Computing $\mathpzc{a}, \mathpzc{b}$ as in \eqref{bconst}, \eqref{a}, the reaction functions $f_j$ read:
$$
f_j(x)=
\begin{cases}
256 \big( \frac{x_1}{1+7x_1}\big)^2 \quad \quad \quad \text{if $j=2$};\\
16\frac{x_i}{1+3x_i} \quad \quad \quad \quad \quad \text{if $j = 2i-1$ or $j = 2i$, $j\neq 2$}.\\
\end{cases}
$$
We compute the Jacobian $G(\mathbf{\bar{r}}')$ at $\mathbf{\bar{r}}'=\mathbf{1}$.
\begin{equation*}
G(\mathbf{1})=
\begin{pmatrix}
-3 & 1 & 0 &...& 0 & 0 & 2\\
1 & -2 & 1 &...& 0 & 0 & 0\\
0 & 1 &  -2 &...& 0 & 0 & 0\\
... & ... & ... & ... & ... & ... & ...\\
0 & 0 & 0 &... & -2 & 1 & 0\\
0 & 0 & 0 &... & 1 & -2 & 1\\
1 & 0 & 0 &... & 0 & 1& -2\\
\end{pmatrix},
\end{equation*}
with right kernel vector $v=(1, 1, ... , 1)^T$ and left kernel vector $w=(M, M+1,..., 2M-1)^T$. Let us first check the condition (SN3),
$$w^T \; \partial^2_x g(x) [v,v] = w^T \sum_i \frac{\partial^2 g}{(\partial x_i)^2} (v_i)^2.$$
A simple computation shows:
\begin{equation*}
(M, M+1, ... , 2M-1)\; \frac{\partial^2 g}{(\partial x_i)^2} =f''_{(2i-1) m_i m_i} -f''_{(2i) m_i m_i},
\end{equation*}
which is nonzero if and only if $i=1$. In fact, note that $f_{2i}\equiv f_{2i-1}$, unless $i=1$. In the case of $i=1$, we have
$$f''_{1 m_1 m_1} -f''_{2 m_1 m_1}=-\frac{13}{8}+\frac{3}{2}=-\frac{1}{8}\neq 0,$$
and thus
$$w^T \; \partial^2_x g(x) [v,v] = w^T \frac{\partial^2 g}{(\partial x_1)^2} \neq 0.$$
Via Lemma \ref{SN2}, or a direct check, we have that
$$\langle w, S^j\rangle  \neq 0,$$
for any reaction $j$: the condition (SN2) is satisfied. In conclusion, the {saddle-node} bifurcation point can be unfolded along $2M$ different parameters $\mathpzc{b_{m_i}^j}$, for $i=1,...,M$, $j=2i$ or $j=2i-1$.\\

\subsection{Example II: glyoxylate cycle vs TCA cycle in \emph{E.coli}} \label{ex2}

The central carbon metabolism is a fundamental metabolic process in living beings. An important part of this process is the \emph{tricarboxylic acid (TCA) cycle}, a cyclic sequence of reactions generating energy in form of ATP. Described for the first time in 1957 by Kornberg and Krebs, the \emph{glyoxylate cycle} is a suggested variation of the TCA cycle. We refer to \cite{kornkrebs57} for more detailed biological explanations. We consider the network structure combining TCA and glyoxylate cycle, as presented in \cite{kornkrebs57}. {Such model does not take in account outflow reactions, which are crucial for a dynamical analysis and indeed abundantly present in dynamical models of metabolism \cite{Chass2002, Ish07}. Thus we further consider outflow reactions as presented in a general model of the Central Carbon Metabolism \cite{Ish07}}. We show that such structure admits a saddle-node bifurcation. The structure is the following:

\begin{center} 
\includegraphics[scale=0.3]{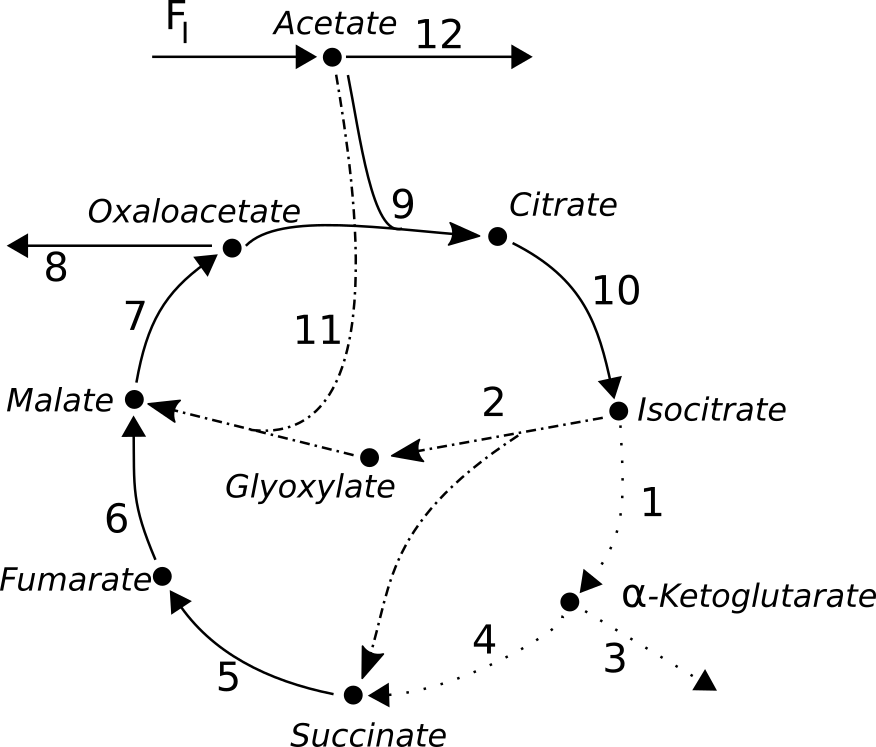}
\end{center}

Above, arrows indicate reactions. Inputs of the arrows are reactants and outputs are products. The continuous arrows refer to reactions present in both TCA and glyoxylate cycle. The sparsely dotted arrows indicate reactions in the TCA cycle not in the glyoxylate cycle: reactions 1, 3, and 4. On the contrary, dotted-dashed arrows indicate reactions of the glyoxylate cycle that do not appear in the TCA cycle: reactions 2 and 11. Reactions 3, 8, and 12 are outflow reactions considered in \cite{Ish07}. The system of differential equations is the following:
$$
\begin{cases}
\dot{x}_A=-r_1(x_A)-r_2(x_A)+r_{10}(x_G);\\
\dot{x}_B=r_1(x_A)-r_3(x_B)-r_4(x_B);\\
\dot{x}_C=r_4(x_B)+r_2(x_A)-r_5(x_C);\\
\dot{x}_D=r_5(x_C)-r_6(x_D);\\
\dot{x}_E=r_6(x_D)-r_7(x_E)+r_{11}(x_H, x_I);\\
\dot{x}_F=r_7(x_E)-r_8(x_F)-r_9(x_F,x_I);\\
\dot{x}_G=r_9(x_F,x_I)-r_{10}(x_G);\\
\dot{x}_H=r_2(x_A)-r_{11}(x_H, x_I);\\
\dot{x}_I=-r_9(x_F,x_I)-r_{11}(x_H, x_I)-r_{12}(x_I)+F_I,
\end{cases}\text{with}\begin{pmatrix}
A\\B\\C\\D\\E\\F\\G\\H\\I
\end{pmatrix}=
\begin{pmatrix}
{\footnotesize\text{\emph{Isocitrate}}}\\ 
{\footnotesize \text{\emph{$\alpha$-Ketoglutarate}}}\\
{\footnotesize\text{\emph{Succinate}}}\\ 
{\footnotesize\text{\emph{Fumarate}}}\\
{\footnotesize\text{\emph{Malate}}}\\
{\footnotesize\text{\emph{Oxaloacetate}}}\\
{\footnotesize\text{\emph{Citrate}}}\\ 
{\footnotesize\text{\emph{Glyoxylate}}}\\
{\footnotesize\text{\emph{Acetate}}}
\end{pmatrix}.
$$
Fix arbitrarily $\bar{r}_3, \bar{r}_4, \bar{r}_8, \bar{r}_{12} > 0$. The equilibrium constraints are:\\
\begin{minipage}{.5\textwidth}
\begin{equation}\label{eqcyc}
\begin{cases}
F_I=3\bar{r}_3+\bar{r}_4+2\bar{r}_8+\bar{r}_{12};\\ \bar{r}_1=\bar{r}_3+\bar{r}_4;\\ \bar{r}_2=\bar{r}_3+\bar{r}_8;\\ \bar{r}_3=\bar{r}_3;\\ \bar{r}_4=\bar{r}_4;\\ \bar{r}_5=\bar{r}_3+\bar{r}_4+\bar{r}_8;\\ \bar{r}_6=\bar{r}_3+\bar{r}_4+\bar{r}_8;\\
\end{cases}
\end{equation}
\end{minipage}%
\begin{minipage}{.5\textwidth}
\begin{equation*}
\begin{cases}
\bar{r}_7=2\bar{r}_3+\bar{r}_4+2\bar{r}_8;\\ \bar{r}_8=\bar{r}_8;\\ \bar{r}_9=2\bar{r}_3+\bar{r}_4+\bar{r}_8;\\ \bar{r}_{10}=2\bar{r}_3+\bar{r}_4+\bar{r}_8;\\ \bar{r}_{11}= \bar{r}_3+\bar{r}_8;\\ \bar{r}_{12}=\bar{r}_{12}.\\
\end{cases}
\end{equation*}
\end{minipage}\\

We identify an SN-pair of Child Selections. Consider
$$\mathbf{J}_1(A,B,C,D,E,F,G,H,I)=(1,3,5,6,7,9,10,11,12),$$
and
$$\mathbf{J}_2(A,B,C,D,E,F,G,H,I)=(2,3,5,6,7,9,10,11,12).$$
$\mathbf{J}_1$ and $\mathbf{J}_2$ are
\begin{enumerate}
\item at minimal distance 1: only the species $A$ is such that $\mathbf{J}_1(A)\neq \mathbf{J}_2(A)$;
\item $\alpha_{\mathbf{J}_1}=-1$ and $\alpha_{\mathbf{J}_2}=+1$;
\item the Partial Child Selection 
$$\mathbf{J}^{\vee A}(B,C,D,E,F,G,H,I)=(3,5,6,7,9,10,11,12)$$
has nonzero coefficient 
$$\beta_{\mathbf{J}^{\vee A}}=1\neq 0.$$
\end{enumerate}
Hence $\mathbf{J}_1$ and $\mathbf{J}_2$ form an SN-pair of Child Selections. The system thus admits a saddle-node bifurcation according to a parametrization of either reaction $1$ or $2$. Reactions $1$ and $2$ mark the difference between the TCA cycle and the glyoxylate cycle. Even when endowed with Michaelis-Menten kinetics, via Theorem \ref{SNMM}, the system exhibits a saddle-node bifurcation. In this case, it is required to choose $\bar{r}_1 \neq \bar{r}_2$, which is allowed by the equilibrium constraints \eqref{eqcyc}.\\

\subsection{Example III: Mass action}\label{ex3}

The law of mass action for a reaction $j$ assumes
\begin{equation}\label{ma}
f_j(x):=k_j \prod_{m\in\mathbf{M}}  x_m^{s^j_m},
\end{equation}
where $k_j >0$ is a positive constant and $s^j_m$ is the stoichiometric coefficient of the species $m$ as reactant of the reaction $j$. Assuming \eqref{ma} for all reaction functions $f_j$ translates \eqref{maineq} into a polynomial system. In contrast to Michaelis-Menten and Hill, a striking feature is that each reaction function $f_j$ is parametrized by only one parameter $k_j$. This impedes our approach, which assumes enough parametric freedom to discuss separately the equilibrium constraints from the bifurcation constraints: it is not the case for mass action. In particular, there is not enough choice of parameters to harness the value of the derivative 
\begin{equation} \label{mader}
f'_{jm}(x)=\;{s^j_m}\; x_m^{(s^j_m - 1)} \; k_j \prod_{n \neq m}  x_n^{s^j_n}=\frac{s^j_m}{x_m} f_j,
\end{equation}
once the value $r_j=f_j(x)$ is fixed. It is not possible using without concern the value $x_m$ as a parameter, as $x_m$ appears also in the mathematical expression of any other derivative $f_{hm}$, for a reaction $h \neq j$ where $m$ participates as reactant. For this reason, the mass-action case deserves further dedication and work, not addressed in this paper. Nevertheless, we derive two observations and produce an example of a network undergoing a saddle-node bifurcation under the assumption of mass action. We keep this discussion as informal and self-contained as possible.\\ 

Essentially, our scheme to detect saddle-node bifurcations is to find two Child Selections $\mathbf{J}_1$, $\mathbf{J}_2$ with opposite behavior and at minimal distance. Firstly, a trivial necessary condition to have a saddle-node is the nonlinearity of the system: for linear systems, condition (SN3) of Theorem \ref{SN} is never satisfied. Under the assumption of mass action, this requires reactions with more than one reactant or with a stoichiometric coefficient bigger than one for the unique reactant. In spirit with the present results, the nonlinearity should be precisely in reactions $j$ with $\mathbf{J}_1(m^*)=j \neq \mathbf{J}_2(m^*)$, for a species $m$. Secondly, in the mass-action case, we are not able to discuss independently equilibrium constraints and bifurcation constraints, as shown in \eqref{mader}. Thus, it is helpful having constant inflow reactions to the species $m^*$ with $\mathbf{J}_1(m^*)=j \neq \mathbf{J}_2(m^*)$. The constant inflow reactions help solve the equilibrium equations, but they do not play any role in the bifurcation conditions, disappearing upon differentiation. 
Following these two observations we present the following example:
$$\underset{F_A}{\longrightarrow} \quad A \quad \underset{1}{\longrightarrow}$$
$$\underset{F_B}{\longrightarrow} \quad B \quad \underset{3}{\longrightarrow}$$
$$\quad A+B  \quad  \underset{2}{\longrightarrow}\quad C \quad \underset{4}{\longrightarrow} \quad 2A+2B$$
where reactions $F_A, F_B$ are inflows and reactions $1, 3$ are outflows to species $A, B$, respectively. From the network we derive the following system of differential equations:
\begin{equation*}
{\small
\begin{cases}
\dot{x}_A=-f_1(x_A)-f_2(x_A,x_B)+2f_4(x_C)+F_A \;=-k_1 x_A - k_2 x_A x_B + 2 k_4 x_C + F_A;\\
\dot{x}_B= -f_3(x_B)-f_2(x_A x_B) +2f_4(x_C)+F_B \;\;=-k_3 x_B - k_2 x_A x_B + 2 k_4 x_C + F_B;\\
\dot{x}_C= f_2(x_A x_B)-f_4(x_C)\quad\quad\quad\quad\quad \quad\quad\quad\;\;=k_2 x_A x_B - k_4 x_C.\\
\end{cases}}
\end{equation*}
Note that reaction $2$ has two reactants, and hence the function $f_2=k_2 x_A x_B$ is nonlinear. Fix arbitrarily $\bar{r}_1, \bar{r}_2, \bar{r_3}>0$ such that $\bar{r}_2<\bar{r}_1, \bar{r}_3$; the equilibrium constraints read:
$$
\begin{pmatrix}
F_A\\
F_B\\
f_1\\
f_2\\
f_3\\
f_4\\
\end{pmatrix}=
\begin{pmatrix}
\bar{r}_1-\bar{r}_2\\
\bar{r}_3-\bar{r}_2\\
\bar{r}_1\\
\bar{r}_2\\
\bar{r}_3\\
\bar{r}_2\\
\end{pmatrix}
\quad$$
The Jacobian of the system is 
$$
G=
\begin{pmatrix}
-k_1-k_2x_B & -k_2x_A & 2k_4 \\
-k_2 x_B & -k_3 -k_2x_A& 2 k_4\\
k_2x_B & k_2x_A& - k_4\\
\end{pmatrix},$$ 
with $\operatorname{det}G=-k_1k_3k_4+k_1k_2 k_4 x_A + k_2 k_3 k_4 x_B$. Let us consider the point $\bar{x}=(\bar{x}_A, \bar{x}_B, \bar{x}_C)=(1,1,1)$, and the rates $k_1=k_3=2$, $k_2, k_4=1$. This solves $\operatorname{det}G = 0$ with a simple eigenvalue zero and fixes $F_A=F_B=1$. Let us consider $k_2$ as a bifurcation parameter, and thus $k^*_2=1$ as its bifurcation value. The Jacobian G at the bifurcation point reads:
$$
G|_{(\bar{x},k_2^*)}=
\begin{pmatrix}
-3& -1 & 2 \\
-1 & -3& 2\\
1 & 1 & - 1\\
\end{pmatrix}
$$
with right kernel vector $v=(1,1,2)^T$ and left kernel vector $w=(1,1,4)$. Condition (SN2) is satisfied:
$$\langle w, \partial_{k_2} g \rangle =(1,1,4) (-1,-1,1)^T \neq 0,$$
as well as (SN3):
$$ w^T \; \partial^2_x g[v,v]= (1,1,4) (-2, -2, 2)^T \neq 0.$$ 
The fact that up to a constant the two conditions (SN2) and (SN3) are the same is not a coincidence, but the central idea of this example: the only nonlinear reaction function $f_2$ is also the only reaction function where the bifurcation parameter appears, hence $\partial_{k_2} g$ must be parallel to $\partial^2_x g[v,v]$. In conclusion, for $k_1=k_3=2$ and ${F_A=F_B=k_4=1}$, the equilibrium $(x_A, x_B, x_C)=(1,1,1)$ undergoes a saddle-node bifurcation for the parameter $k_2=1$. In particular, for $k_2 < 1$ we have multistationarity.

\subsection{Example IV: Degenerate saddle-node for Michaelis-Menten kinetics} \label{ex5}

This example presents a network that, when endowed with Michaelis-Menten kinetics, admits an equilibrium $\bar{x}$ satisfying conditions (SN1) and (SN2) but for which (SN3) is never satisfied. The network is the following:
$$\underset{0}{\longleftarrow} \quad A \quad\underset{1}{\longrightarrow}\quad B \quad \underset{2}{\longrightarrow} \quad 2A$$
where reaction $0$ is an outflow from $A$. The system possesses only two Child Selections $\mathbf{J}_1$ and $\mathbf{J}_2$: 
$$\mathbf{J}_1(A,B)=(0,2) \quad\quad\quad\text{and}\quad\quad\quad \mathbf{J}_2(A,B)=(1,2).$$
$\mathbf{J}_1$ and $\mathbf{J}_2$ form an SN-pair of Child Selections. In fact, they are at minimal distance $d(\mathbf{J}_1, \mathbf{J}_2)=1$, $\alpha_{\mathbf{J}_1} \alpha_{\mathbf{J}_2} = 1 \cdot (-1) = -1$, and the unique Partial Child Selection $\mathbf{J}^{\vee A}(B)=2$ is nonzero. Theorem \ref{SNmy} guarantees that there exist monotone chemical functions such that the associated dynamical system undergoes a saddle-node bifurcation according to a parametrization of the reaction $\eta=0$ or $\eta=1$. However, when restricted to Michaelis-Menten kinetics the saddle-node is always degenerate. Consider the associated system
\begin{equation*}
\begin{cases}
\dot{x}_A= -r_0(x_A)-r_1(x_A)+2r_2(x_B);\\
\dot{x}_B= r_1(x_A)-r_2(x_B).\\
\end{cases}
\end{equation*}
The equilibrium constraints fix 
$\bar{r}_0=\bar{r}_1=\bar{r}_2,$ and hence condition \eqref{MMcond} is never satisfied.\\

However, the degenerate situation is easily fixable: let us consider the same system with an added inflow to species $A$:
$$\quad \underset{F_A}{\longrightarrow} \quad A.$$
The system of ODEs now reads
\begin{equation*}
\begin{cases}
\dot{x}_A= -r_0(x_A)-r_1(x_A)+2r_2(x_B)+F_A;\\
\dot{x}_B= r_1(x_A)-r_2(x_B),\\
\end{cases}
\end{equation*}
with the equilibrium constraints:
$$\bar{r}_1=\bar{r}_2=\bar{r}_0-F_A.$$
Thus, for $0 < F_A < \bar{r}_0$, we have now equilibria for which $\bar{r}_0 \neq \bar{r}_1$. The case $F_A = 0$ recovers the degenerate example. If $\bar{r}_0 \neq \bar{r}_1$ we have a nondegenerate saddle-node bifurcation, even in the Michaelis-Menten case.

\section{Discussion}\label{discussion}

We have presented a comprehensive saddle-node bifurcation analysis for chemical reaction networks. Via a symbolic approach, we have analyzed which networks admit the occurrence of a bifurcation behavior. This work has two direct consequences. Theoretically, we have described the structures that guarantee that the network can sustain multistationarity. Practically, we have identified the proper parameters to unfold a saddle-node bifurcation.\\

\begin{figure}
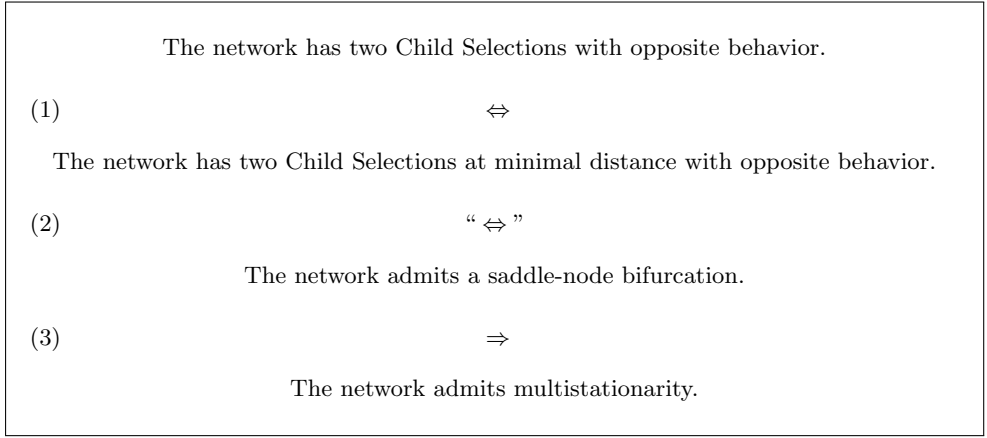

\begin{framed}
{\small
$$\text{The network has two Child Selections with opposite behavior.}$$
\begin{equation}\tag{1}
\Leftrightarrow
\end{equation}
$$\text{The network has two Child Selections at minimal distance with opposite behavior.}$$
\begin{equation}\tag{2}
``\Leftrightarrow"
\end{equation}
$$\text{The network admits a saddle-node bifurcation.}$$
\begin{equation}\tag{3}
\Rightarrow
\end{equation}
$$\text{The network admits multistationarity.}$$}
\end{framed}
\caption{The conceptual map of this paper. Implication (1) is Lemma \ref{mainlemma}. Implication (2) is technically only a necessary condition $\Leftarrow$, but many ``pathological" counterexamples are not biologically relevant, hence the implication is expected to be also sufficient in realistic networks. Implication (3) concludes the logical chain by assessing multistationarity.}
\label{figure}
\end{figure}

The key structure we have described are \emph{SN-pairs} of Child Selections, i.e., two nonzero Child Selections $(\mathbf{J}_1, \mathbf{J}_2)$ that satisfy three conditions:
\begin{enumerate}
\item they are at minimal distance $\delta$;
\item their behavior coefficient is opposite in sign: $\alpha_{\mathbf{J}_1}\alpha_{\mathbf{J}_2}<0$;
\item a technical condition excluding multiple eigenvalues zero.
\end{enumerate}
If the network possesses an SN-pair of Child Selections, then a saddle-node bifurcation occurs for a choice of monotone chemical functions $f(x,\lambda)$ parametrized by a single parameter $\lambda$. The bifurcation parameter $\lambda$ parametrizes only the function $f_\eta$ of an arbitrary reaction $\eta=\mathbf{J}_1(m^*)\neq \mathbf{J}_2(m^*)$, for one of the $\delta$ species $m^*$ with $\mathbf{J}_1(m^*)\neq \mathbf{J}_2(m^*)$. The existence of a pair of Child Selections satisfying conditions (1) and (2) above is only necessary for a bifurcation behavior. However, the identification of a counterexample \cite{V22} to the sufficiency of conditions (1) and (2) suggests that this is an issue for mathematicians, with seemingly no biological relevance. That is, for realistic biological networks, conditions (1) and (2) above essentially characterize the bifurcation behavior. On the other hand, the three conditions (1)--(3) are technically only sufficient. See Figure \ref{figure} for a conceptual map.\\

Our symbolic approach considers the derivatives $\mathbf{r}'$ of the reaction functions as positive independent variables. We have applied a geometrical perturbation argument with an $\varepsilon$-rescaling of the variables $\mathbf{r}'$. At the limit $\varepsilon=0$, the only nonzero variables are the ones identified by the SN-pair $(\mathbf{J}_1, \mathbf{J}_2)$ of Child Selections, that is:
$$r'_{jm} \neq 0 \quad\quad\quad \Leftrightarrow \quad\quad\quad (j,m)=(\mathbf{J}_i(m),m)\quad \text{ for $i=1,2$.}$$
In this sense, the bifurcation behavior of the SN-pair of Child Selections is \emph{inherited} by the full network. Inheritance of dynamical features figures in recent works by Banaji \cite{Ba18, BaPa18}. However, these works have not yet discussed the inheritance of bifurcation behavior, and focused on modified networks, rather than identifying some leading subnetworks that encode the dynamics, as we did. In particular, any network for which a pair of Child Selections is an SN-pair admits a saddle-node bifurcation. Corollaries \ref{mindist1} exploits and clarifies this idea, and can be used to find \emph{small and simple network motifs} for saddle-node bifurcation in larger networks.\\

Network motifs connected to multistationarity have been discussed in the literature on various levels of abstraction and empiricism. A central role is often claimed by \emph{autocatalysis} \cite{Schu19}. Autocatalytic reactions are those in which at least one of the products is also a reactant. For example, the reversible reaction
\begin{equation*}\label{jaut}
m \quad \overset{j^+_{aut}}{\underset{j^-_{aut}}{\myrightleftarrows{\rule{0.5cm}{0cm}}}} \quad 2m
\end{equation*}
is autocatalytic. We refer to reaction $j^+_{aut}$ as \emph{positively autocatalytic} and to $j^-_{aut}$ as \emph{negatively autocatalytic}. These ideas can be generalized to autocatalytic sequences of reactions (feedback loops) or autocatalytic networks, leading to a more general concept of autocatalysis, which is still under debate \cite{DefAut20}. We do not enter here such a formal discussion, but we observe that the results are consistent with the many independent observations on the centrality of autocatalysis.  In particular, positive autocatalysis can trigger a saddle-node bifurcation and consequent multistationarity. Example \ref{ex1} shows how reversible feedback loops admit saddle-node bifurcations. Such feedback loops are a simple generalization of autocatalytic reactions and their connection with multistationarity has been pioneered by Thomas \cite{Thom81} in a related context. Example \ref{ex2} identifies a saddle-node bifurcation in the central carbon metabolism of \emph{E. coli}. Along the glyoxylate cycle, one molecule of \emph{Isocitrate} transforms into two molecules, similarly to $j^+_{aut}$. Example \ref{ex3} contains two reactions
$$\quad A+B  \quad  \underset{2}{\longrightarrow}\quad C \quad \underset{4}{\longrightarrow} \quad 2A+2B,$$
with a clear analogy to $j^+_{aut}$. Example \ref{ex5} contains a positive feedback loop of two reactions
$$A\quad\underset{1}{\longrightarrow}\quad B \quad \underset{2}{\longrightarrow} \quad 2A.$$
More simply, let us consider a toy network $\mathbf{\Gamma}$ with a single species $A$. The behavior coefficient of any Child Selection on $\mathbf{\Gamma}$ is then just 
$$\alpha_\mathbf{J}=\operatorname{det}S^\mathbf{J}=S_{A \mathbf{J}(A)}=\tilde{s}^{\mathbf{J}(A)}_A- s^{\mathbf{J}(A)}_A.$$
Clearly, a Child Selection $\mathbf{J}$ is bad if and only if $\mathbf{J}(A)$ is positively autocatalytic. This observation can be generalized to any dimension and highlights the connection between positive autocatalysis and bad Child Selections. Previous work \cite{VGB20} argued on the prevalence of good Child Selections in biochemical networks. In light of such an argument, the mere presence of positive autocatalysis points to the first step of the logical chain of Figure \ref{figure} and the consequent possibility of multistationarity. The formalization of these presented arguments will be included in future work.\\

Realistic kinetic models of biochemical networks typically comprise different types of kinetics \cite{Chass2002}. In this paper, we have explicitly discussed the case of two kinetics of interest: Michaelis-Menten and Hill kinetics, and identified the proper bifurcation parameters.  In the presence of a SN-pair of Child Selections $(\mathbf{J}_1, \mathbf{J}_2)$, the bifurcation parameter is $b^\eta_{m^*}$ in \eqref{MM} for any reaction $\eta$ and species $m^*$ such that $\eta=\mathbf{J}_1(m^*)\neq \mathbf{J}_2(m^*)$. We have shown that the parametric richness of Hill kinetics always guarantees a saddle-node bifurcation behavior. On the contrary, for Michaelis-Menten kinetics we can only guarantee the existence of an equilibrium such that conditions (SN1)-(SN2) of Theorem \ref{SN} are satisfied, not necessarily (SN3). However, Theorem \ref{SNMM} provides a sufficient condition \eqref{MMcond} for a nondegenerate saddle-node bifurcation in the case of distance $d(\mathbf{J}_1, \mathbf{J}_2)=1$. Condition \eqref{MMcond} shows that the degeneracy is quite unlikely, as it requires a uniquely determined equilibrium constraint \eqref{eqconst}. Example \ref{ex5} shows how an inflow reaction to $m^*$ already guarantees that there exists a choice of equilibrium fluxes $\mathbf{r}$ such that \eqref{MMcond} is satisfied. Finally, the analysis for mass action kinetics is at present inconclusive, even though Example \ref{ex3} presents a mass-action system undergoing a saddle-node bifurcation according to the same structural intuition of the present paper. Further work is needed to clarify it, along the lines of the present contribution.\\

A natural combinatorial question arises with regard to SN-pairs of Child Selections $(\mathbf{J}_1,\mathbf{J}_2)$: which structure is required for the case of minimal distance $d(\mathbf{J}_1,\mathbf{J}_2)>1$? In Example \ref{ex1} we have presented the case of reversible feedback cycles of length $M$, which possesses an SN-pair of Child Selections at distance $M$. In a work in preparation, we will characterize SN-pairs of Child Selections in terms of a reciprocal permutation structure, generalizing the concept of reversible cycles. This observation stresses how saddle-node bifurcations are triggered either by SN-pairs at distance $d=1$, or by special (thus recognizable!) structures. Note that even for Michaelis-Menten kinetics we have a sufficient condition for saddle-node for the most relevant case $d=1$.\\

In conclusion we make a non-mathematical consideration. It is a strong impression of the author that the bifurcation behavior is essentially characterized in realistic biological networks by simple and recognizable structures, even though mathematics includes much more complex options, and requires a thorough analysis and exclusion of pathological cases. Such pathological cases are of minimal interest for biologists but enhance the technical difficulties of the proofs and diminish the verbal strength of the results. This empirical consideration calls for stronger interaction between mathematicians and theoretical biologists to synthesize the results, with a skimming of biologically irrelevant (but nevertheless mathematically challenging and intriguing) cases, with the goal of obtaining a clearer picture of the bifurcation behavior in real-world biochemical networks.

\section{Proofs}\label{proofs} \textcolor{white}{,}\\

\textbf{Proofs of Section 5 and Theorem \ref{det0}}: before proving Lemma \ref{mindistsol}, we prove a further lemma.
\begin{lemma}\label{Lemma2}
Let $\mathbf{J}_1, \mathbf{J}_2$ be two Child Selections at minimal distance $\delta$. Then, for every other nonzero Child Selection $\mathbf{J}_3$, there exists a species $m^*$ such that
\begin{equation}\label{thesis}
\begin{cases}
\mathbf{J}_3(m^*)\neq \mathbf{J}_1(m^*);\\
\mathbf{J}_3(m^*)\neq \mathbf{J}_2(m^*).\\
\end{cases}
\end{equation}
\end{lemma}

\proof
Let $\mathbf{D} \subseteq \mathbf{M}$ be the set of species $m$ such that $\mathbf{J}_1(m) \neq \mathbf{J}_2(m)$. If there exists $m^* \in \mathbf{M} \setminus \mathbf{D}$ such that $\mathbf{J}_1(m^*) \neq \mathbf{J}_3(m^*)$
we are done, since together with $\mathbf{J}_1(m^*)=\mathbf{J}_2(m^*)$ it implies \eqref{thesis}. Consequently, let us assume $\mathbf{J}_1(m)=\mathbf{J}_2(m) = \mathbf{J}_3(m)$, for every $m \in \mathbf{M} \setminus \mathbf{D}$. Take any $m^* \in \mathbf{D}$. We have \eqref{thesis}. Indeed, assume \eqref{thesis} does not hold. Without loss of generality we have $\mathbf{J}_3(m^*)= \mathbf{J}_1(m^*),$ which implies
$$d(\mathbf{J}_1, \mathbf{J}_3) \le \delta - 1 < \delta,$$
and contradicts the assumption of $\delta$ being the minimal distance.
\endproof

\proof[Proof of Lemma \ref{mindistsol}]
For any $r'_{jm} \neq r'_{\mathbf{J}_i(m)m}$, for $i=1, 2$, fix the value
$$r'_{jm}=\varepsilon \; \bar{r}_{jm},$$
where $\bar{r}_{jm} > 0$ is any positive value. Via Lemma \ref{Lemma2}, $P(\mathbf{r}')$ now takes the form:
\begin{equation}\label{poleps}
P(\mathbf{r}')=\alpha_{\mathbf{J}_1} \mathbf{r}'[\mathbf{J}_1] + \alpha_{\mathbf{J}_2} \mathbf{r}'[\mathbf{J}_2] + q(\varepsilon),
\end{equation}
where $q(\varepsilon)$ indicates all the summands where $\varepsilon$ appears at least linearly. Lemma \ref{Lemma2} guarantees indeed that there are no other nonzero summands. In particular, $q(\varepsilon)$ is a polynomial in $\varepsilon$ with zero constant term, i.e. $q(0)=0$. 
\endproof

\proof[Proof of Lemma \ref{mainlemma}]
The implication $\Leftarrow$ is trivial. We prove the implication $\Rightarrow$. Consider the following two sets of Child Selections:
$$\mathpzc{G} = \{ \mathbf{J} \; | \; \alpha_\mathbf{J}=(-1)^M \} \quad \text{and} \quad \mathpzc{B} = \{ \mathbf{J} \;  |\; \alpha_\mathbf{J}=(-1)^{M-1} \}.$$
Since the total number of Child Selections is finite, $\mathpzc{G}$ and $\mathpzc{B}$ are finite sets and, by assumption, both are nonempty. We can define the distance $d$ of the two sets as:
$$d(\mathpzc{G},\mathpzc{B}):=\operatorname{inf}_{\mathbf{J}_1 \in \mathpzc{G}, \mathbf{J}_2 \in \mathpzc{B}} d(\mathbf{J}_1, \mathbf{J}_2)=\operatorname{min}_{\mathbf{J}_1 \in G, \mathbf{J}_2 \in B} d(\mathbf{J}_1, \mathbf{J}_2)=d(\mathbf{J}_3, \mathbf{J}_4),$$
for some Child Selections $\mathbf{J}_3$ and $\mathbf{J}_4$ with $\alpha_{\mathbf{J}_3}\alpha_{\mathbf{J}_4} < 0.$
The Child Selections $\mathbf{J}_3$ and $\mathbf{J}_4$ are at minimal distance, by construction.
\endproof

\proof[Proof of Theorem \ref{det0}] The implication $\Rightarrow$ is trivial. We prove the implication $\Leftarrow$. We apply Lemma \ref{mainlemma} to find two Child Selections $\mathbf{J}_3, \mathbf{J}_4$ at minimal distance with $\alpha_{\mathbf{J}_3}\alpha_{\mathbf{J}_4}< 0$. We recall the set $\mathbf{D} \subseteq \mathbf{M}$ of species $m$ such that $\mathbf{J}_3(m) \neq \mathbf{J}_4(m)$. The cardinality of $\mathbf{D}$ is $\delta=d(\mathbf{J}_3,\mathbf{J}_4)$. Without loss of generalities let us consider $\mathbf{D}=\{m_1, ..., m_\delta\}$. \\

We consider the $\varepsilon$-rescaling of Lemma \ref{mindistsol}. We define $m^*:=m_1$, $\eta:= \mathbf{J}_3(m^*)$ and $\rho:=r'_{\eta m^*}$. We want to solve $P(\mathbf{r}')=0$ with respect to the variable $\rho$. To this goal, for any $r'_{jm} = r'_{\mathbf{J}_i(m)m}$, for $i=3, 4$, $r'_{jm}\neq \rho$, we fix the value
$$r'_{jm} := \bar{r}'_{jm},$$
where $\bar{r}'_{jm} > 0$ is any positive value. Now the polynomial \eqref{poleps} reads as a bivariate polynomial in the two variables $\rho$ and $\varepsilon$:
$$P(\mathbf{r}')=P(\rho,\varepsilon)=\alpha_{\mathbf{J}_3} \rho \; \bar{\mathbf{r}}'[\mathbf{J}_3 \setminus \eta] + \alpha_{\mathbf{J}_4} \bar{\mathbf{r}}'[\mathbf{J}_4] + q(\varepsilon),$$
where $\bar{\mathbf{r}}'[\mathbf{J}_3 \setminus \eta]$ indicates the monomial $\bar{\mathbf{r}}'[\mathbf{J}_3 \setminus \eta]:=\frac{\bar{\mathbf{r}}'[\mathbf{J}_3]}{\rho}.$ The value 
$$\rho^*:= -\frac{\alpha_{\mathbf{J}_4}\bar{\mathbf{r}}'[\mathbf{J}_4]}{ \alpha_{\mathbf{J}_3}\bar{\mathbf{r}}'[\mathbf{J}_3 \setminus\eta]},$$
is positive by assumption, since $\alpha_{\mathbf{J}_3}\alpha_{\mathbf{J}_4}<0$. Moreover,
$$P(\rho^*, 0)=0.$$
We apply the Implicit Function Theorem to show that a positive solution $\rho^*(\varepsilon)$ persists also in a positive neighborhood $\varepsilon>0$. We check
$$P_\rho (\rho, \varepsilon)|_{(\rho^*,0)} \neq 0.$$
Indeed,
$$P_\rho (\rho, \varepsilon)|_{(\rho^*,0)} = \alpha_{\mathbf{J}_3}\bar{\mathbf{r}}'[\mathbf{J}_3 \setminus \eta] + q_\rho(\varepsilon)|_{(\rho^*,0)}=\alpha_{\mathbf{J}_3}\bar{\mathbf{r}}'[\mathbf{J}_3 \setminus\eta] \neq 0.$$
Hence, there exists a positive solution $\rho^*(\varepsilon)$ for $\varepsilon > 0$.
\endproof
\begin{remark}
In the proof of Theorem \ref{det0}, the choice of $\rho$ is arbitrary. We could argue analogously, by picking any $\rho^*:=r'_{\mathbf{J}_i(m)m}$,  $i=3,4$, for any $m$ such that $\mathbf{J}_3(m)\neq\mathbf{J}_4(m)$.
\end{remark}
\begin{remark}
Theorem \ref{det0} can also be proved via intermediate value theorem, in simpler setting. However, the perturbation construction in the presented proof is central to the development of the following results of this paper.
\end{remark}

\textbf{Proofs of Section 6}

\proof[Proof of Lemma \ref{simple0}]

Via condition (1) of Definition \ref{snpairs} we consider again the $\varepsilon$-rescaling of Lemma \ref{mindistsol} and proceed as in the proof of Theorem \ref{det0}. Without loss of generality, assume again that the $\delta$ species $m$ such that $\mathbf{J}_1(m) \neq  \mathbf{J}_2(m)$ are $m_1, ... , m_\delta$, and that $\tilde{m}=m_1$ in the Definition \ref{snpairs}. Define again $\eta:= \mathbf{J}_1(\tilde{m})$, and $\rho:=r'_{\mathbf{J}_1(\tilde{m})\tilde{m}}$. With abuse of notation, let $\mathpzc{r'}$ indicate all the variables 
$$\mathpzc{r}':=\{r'_{jm}\text{ such that } r'_{\mathbf{J}_i(m)m},\; i=1,2, \;r'_{jm}\neq \rho \}.$$ 
This way the polynomial $P(\mathbf{r}')$ reads as $P(\rho, \mathpzc{r'}, \varepsilon)$, in the variables $\rho, \mathpzc{r'}, \varepsilon$. Let now $\rho^*(\mathpzc{r}')$ indicate the function of $\mathpzc{r'}$
$$\rho^*(\mathpzc{r'}):= -\frac{\alpha_{\mathbf{J}_2}{\mathbf{r}}'[\mathbf{J}_2]}{ \alpha_{\mathbf{J}_1}{\mathbf{r}}'[\mathbf{J}_1 \setminus\eta]}.$$
Condition (2) of Definition \ref{snpairs} guarantees positivity of $\rho^*(\mathpzc{r'})>0$, and thus
$$P(\rho^*(\mathpzc{r'}), \mathpzc{r'},0)=0, \quad\text{for any choice of $\mathpzc{r'}>0$}.$$
On the other hand, via \eqref{23}, 
$$
\operatorname{Adj}G^{\tilde{m}}_{\tilde{m}}= \sum \limits_{\mathbf{J}^{\vee \tilde{m}}}\beta_{\mathbf{J}^{\vee \tilde{m}}} \; \mathbf{r}'[\mathbf{J}^{\vee \tilde{m}}],
$$
implies that $\rho$ does not appear in $\operatorname{Adj}G^{\tilde{m}}_{\tilde{m}}$. Moreover, Condition $(3)$ of Definition \ref{snpairs} guarantees the existence of a nonzero Partial Child Selection $\mathbf{J}^{\vee \tilde{m}}$, i.e. $ \beta_{\mathbf{J}^{\vee \tilde{m}}}\neq 0$. Thus,
$$\operatorname{Adj}G^{\tilde{m}}_{\tilde{m}} \not \equiv 0.$$
In conclusion, we have that 
$$A(\rho^*(\mathpzc{r'}), \mathpzc{r'}, 0)=\operatorname{tr} \operatorname{Adj} G (\rho^*, \mathpzc{r'}, 0) \not \equiv 0.$$
Hence, there exists a choice of $\bar{\mathpzc{r}}'$ such that
\begin{equation*}
\begin{cases}
P(\rho^*(\bar{\mathpzc{r}}'), \bar{\mathpzc{r}}', 0) = 0;\\
A(\rho^*(\bar{\mathpzc{r}}'), \bar{\mathpzc{r}}', 0)\neq 0.
\end{cases}.
\end{equation*}
As in the proof of Theorem \ref{det0}, we apply the implicit function theorem obtaining a solution of $P(\rho^*( \bar{\mathpzc{r}}',\varepsilon), \bar{\mathpzc{r}}', \varepsilon)$, for small $\varepsilon$. By continuity, $A(\rho^*, \bar{\mathpzc{r}}', 0)\neq 0$ persists in a neighborhood of $\varepsilon$, implying the existence of $\bar{\mathbf{r}}'$ such that:
\begin{equation*}
\begin{cases}
P(\bar{\mathbf{r}}') = 0;\\
A(\bar{\mathbf{r}}')\neq 0.
\end{cases}
\end{equation*}
\endproof

\textbf{Proofs of Section 7}\\

\nopagebreak

\proof[Proof of Lemma \ref{SN2}]
Let us consider the $\varepsilon$-rescaling of Lemma \ref{mindistsol} for the SN-pair of Child Selections $\mathbf{J}_1$ and $\mathbf{J}_2$. Until the last step of this proof, we consider the limit $\varepsilon=0$. Not to overload notation, we drop the constant reference to $\varepsilon$. At this limit, we recall that the only nonvanishing variables $r'_{jm}$ are the ones such that
$$r'_{jm}=r'_{\mathbf{J}_i(m)m} \quad\quad\quad \text{with $i=1,2$}.$$

Let $\rho:=r'_{\eta m^*}$ and consider the $M \times M$ matrix $G^{\eta}$, whose columns $(G^\eta)^m$ are
\begin{equation*}
(G^\eta)^m=
\begin{cases}
\rho \; S^\eta \quad \quad \quad \quad \;\;\; \text{if $m=m^*$};\\
G^m \quad \quad \quad \quad \quad \text{otherwise}.
\end{cases}
\end{equation*}
We recall that $G^m$ indicates the $m^{th}$ column of $G$ and $S^\eta$ the stoichiometric column corresponding to reaction $\eta$. The polynomial $P^\eta(\mathbf{r}'):= \operatorname{det}G^{\eta}$ is nonzero. In fact
\begin{equation}\label{detmod}
P^\eta(\mathbf{r}')=\sum_{\eta \in \mathbf{J}} \alpha_\mathbf{J} \mathbf{r}'[\mathbf{J}] = \alpha_{\mathbf{J}_1} \mathbf{r}'[\mathbf{J_1}] \neq 0, \quad \text{for any choice of  $\mathbf{r}'[\mathbf{J_1}]>0$}.
\end{equation}
Note that \eqref{detmod} holds only at the limit $\varepsilon=0$. Let again $\mathpzc{r'}$ indicate all the variables 
$$\mathpzc{r}':=\{r'_{jm}\text{ such that } r'_{\mathbf{J}_i(m)m},\; i=1,2, \;r'_{jm}\neq \rho \}.$$ 
As in the proof of Lemma \ref{simple0}, we can choose $\bar{\mathpzc{r}}', \rho^*(\bar{\mathpzc{r}}')$ such that 
$G[(\rho^*(\bar{\mathpzc{r}}'),\bar{\mathpzc{r}}')]$ has an algebraically simple eigenvalue zero. Let $w$ be left kernel vector of $G[(\rho^*(\bar{\mathpzc{r}}'),\bar{\mathpzc{r}}')]$. Note that the $m^{th}$ column of $G$ is the $m^{th}$ column of $G^{\eta}$, except for $m=m^*$. From the nonsingularity of  $G^{\eta}$ we conclude: 
\begin{equation*}
0 \neq w^T G^{\eta}=
\begin{blockarray}{ccccc}
1 & ... & m^* & ... & M \\
\begin{block}{(ccccc)}
  0, & ... & \langle w, \rho S^{\eta}\rangle , & ... & 0  \\
\end{block}
\end{blockarray},
\end{equation*}
implying $\langle w, S^{\eta}\rangle \neq 0$, which persists for small $\varepsilon>0$, by continuity.
\endproof

\proof[Proof of Lemma \ref{sn3noneq0}]

By linearity
\begin{equation*}
\begin{split}
w^T \; \partial^2_x g(\bar{x}, \lambda^*)[v,v] &= w^T \; \sum_{h,k} \frac{\partial^2 g}{\partial x_h \partial x_k}(v_h v_k)\\
&=\sum_{h,k} w^T \; \frac{\partial^2 g}{\partial x_h \partial x_k}(v_h v_k).
\end{split}
\end{equation*}
For any $m$ and $n$, a second derivative $f''_{jmn}$ appears only in the summand $$w^T  \frac{\partial^2 g}{\partial x_m \partial x_n}(v_m v_n).$$ Hence,
$$w^T \; \partial^2_x g(\bar{x}, \lambda^*)[v,v]  \equiv 0 \quad  \quad \Leftrightarrow \quad  \quad \sum_{h,k} w^T \frac{\partial^2 g}{\partial x_m \partial x_n}(v_m v_n)\equiv 0 \text{ for any $m$ and $n$}.$$
We show that there exists $m \in \mathbf{M}$ such that:
\begin{equation}\label{lemmasn3eq0}
w^T \frac{\partial g^2 (\bar{x}, \lambda^*)}{(\partial x_m)^2}(v_m)^2 \not \equiv 0.
\end{equation}
Indeed, since $\bar{\mathbf{r}}'$ is a positive solution of $P(\mathbf{r}')=0$,  there exist $m^*$ and $\eta$ such that
\begin{equation}\label{lemmasn3eq}
\rho^*:=\bar{r}'_{\eta m^*} =  \frac{-\sum_{j \not \in \mathbf{J}} \; \alpha(\mathbf{J}) \; \bar{\mathbf{r}}'[\mathbf{J}]}{\sum_{j \in \mathbf{J}}\; \alpha(\mathbf{J})\; \bar{\mathbf{r}}'[\mathbf{J}\setminus j]},
\end{equation}
with numerator and denominator of the same sign. In the proof of Theorem \ref{det0}, we have shown how to construct such $\rho^*$ from two Child Selections of opposite behavior. Note that the $m^{th}$ column $G^{m}=\partial g/\partial{x_m}$ of the Jacobian $G$ has an identical symbolic structure as the the column vector $g''_{mm}:=\partial^2 g/(\partial{x_m})^2$. Indeed, every first derivative $r'_{jm}$ in $G^m$ is simply substituted with the second derivative $r''_{jmm}$ in $g''_{mm}$. This implies
$$w^T \frac{\partial g^2 (\bar{x}, \lambda^*)}{(\partial x_m)^2} = 0$$
for any $m$, if $r'_{jm}=r''_{jmm},$ for every $m$ and $j$. Let us focus on $m^*$. For $(j,m) \neq (\eta, m^*)$, fix $\bar{r}''_{jmm} = \bar{r}'_{jm}$, and let $\bar{r}''_{\eta m^*} \neq \rho^*$. Clearly, from \eqref{lemmasn3eq} this choice of $\bar{\mathbf{r}}''$ implies
$$w^T \frac{\partial g^2 (\bar{x}, \lambda^*)}{(\partial x_{m^*})^2}(\bar{\mathbf{r}}'') \neq 0.$$
To conclude \eqref{lemmasn3eq0}, we show that $v_{m^*} \neq 0$. Fix $\bar{r}'_{jm}$ for $(j,m) \neq (\eta, m^*)$ as in \eqref{lemmasn3eq} and let $\rho$ be the only variable. Then the univariate polynomial
$$P(\mathbf{r}')=P(\rho)$$
is evaluated zero if and only if $\rho=\rho^*$. But $\rho$ appears only in the $m^*\;^{th}$ column of $G$. Via $G(\mathbf{r}') v= 0$, we have $v_{m^*} \neq 0$.
\endproof

\textbf{Proofs of Section 8} \nopagebreak
\proof[Proof of Theorem \ref{S0MM}]
Via Lemma \ref{simple0}, if the network possesses an SN-pair of Child Selections, then there is a choice of $\bar{x}, \bar{\mathbf{r}}, \bar{\mathbf{r}}'$ such that \eqref{posb1} is satisfied, with the Jacobian $G(\bar{\mathbf{r}}')$ possessing an algebraically simple eigenvalue zero. We can then choose $(\mathbf{a}, \mathbf{b})=(\mathpzc{a},\mathpzc{b})$, with $\mathpzc{a},\mathpzc{b}$ defined as in \eqref{a}, \eqref{bconst}, respectively. The choice of the bifurcation parameter 
$$\lambda:=b^\eta_{m^*} \quad\quad\quad\text{with bifurcation value}\quad\quad\quad \lambda^*:=\mathpzc{b}^\eta_{m^*}$$
implies that 
$$\partial_\lambda g =\frac{\partial f_\eta}{\partial \lambda} S^\eta,$$
where $S^\eta$ is the stoichiometric vector of reaction $\eta$. Since $\frac{\partial f_\eta}{\partial \lambda} \neq 0$, Lemma \ref{SN2} implies that the Michaelis-Menten system possesses a simple eigenvalue zero at $\bar{x}$ for the choice $(\mathbf{a}, \mathbf{b})=(\mathpzc{a},\mathpzc{b})$, which satisfies (SN1) and (SN2) of Theorem \ref{SN}.
\endproof

%Before we go on proving Theorem \ref{SNMM}, let us state and prove a needed Lemma.
%\begin{lemma}
%Let $\mathbf{\Gamma}$ be a network possessing \emph{only} two Child Selections $\mathbf{J}_1$ and $\mathbf{J}_2$ at distance $d=1$. Assume $\alpha_{\mathbf{J}_1}\alpha_{\mathbf{J}_2}<0$. Let $m^*$ be the unique species such that  $\mathbf{J}_1(m^*)=j_1\neq j_2=\mathbf{J}_2$. Let $\bar{x}$ be an equilbrium, and assume $\alpha_{\mathbf{J}_1}f'_{j_1m^*}(\bar{x})=\alpha_{\mathbf{J}_2}f'_{j_2 m^*}(\bar{x})$. Then the Jacobian $G$ is singular with $v$ and $w$ respectively left and right eigenvector and it holds $$w^T \; \partial^2_x g''[v, v] = \alpha_{\mathbf{J}_1}f''_{j_1m^*m^*}(\bar{x})-\alpha_{\mathbf{J}_2}f''_{j_2 m^*m^*}(\bar{x})$$
%\end{lemma}

\proof[Proof of Theorem \ref{SNMM}]
Let us consider the $\varepsilon$-rescaling of Lemma \ref{mindistsol} for the SN-pair of Child Selections $\mathbf{J}_1, \mathbf{J}_2$. Until the last step of this proof, we consider the limit $\varepsilon=0$. Not to overload notation, we drop the constant reference to $\varepsilon$. At this limit, we recall that the only nonvanishing variables $r'_{jm}$ are the ones such that
$$r'_{jm}=r'_{\mathbf{J}_i(m)m} \quad\quad\quad \text{with $i=1,2$}.$$
Note that distance $d(\mathbf{J}_1, \mathbf{J}_2)=1$ and injectivity of Child Selections imply
\begin{equation}\label{mix0}
r'_{j m}r'_{j n}=0
\end{equation}
for any reaction $j$ and $m\neq n$.\\

At $\varepsilon=0$, we can proceed as in the proof of Lemma \ref{simple0}. For a proper choice of $\bar{\mathbf{r}}'$ we obtain an algebraically simple eigenvalue zero. Note that
\begin{equation}\label{detresc0}
P(\mathbf{r}')=\operatorname{det}G(\mathbf{r}')=(\alpha_\mathbf{J_1}r'_{\eta m^*}+\alpha_{\mathbf{J}_2} r'_{j_2 m^*})\mathbf{r}'[\mathbf{J}_1 \setminus \eta]= 0,
\end{equation}
if and only if 
\begin{equation}\label{conddet0}
\alpha_{\mathbf{J}_1} r'_{\eta m^*} =-\alpha_{\mathbf{J}_2} r'_{j_2m^*}.
\end{equation}
Let $v, w$ be right and left kernel vectors of $G(\bar{\mathbf{r}}')$, respectively. Let $G^m$ indicate the $m^{th}$ column of the Jacobian $G$. At the limit $\varepsilon=0$ we have 
\begin{equation*}
\begin{cases}
 G^m = S^{\mathbf{J}_1(m)}  r'_{\mathbf{J}_1(m)m} \quad \quad \quad \text{for $m\neq m^*$};\\
 G^{m^*}= S^{\eta} r'_{\eta m^*} + S^{j_2} r'_{j_2 m^*}.\\
\end{cases}
\end{equation*}
Consequently, $w^T G(\mathbf{r}') = 0$ yields
\begin{equation*}
\begin{cases}
w^T G^m = w^T S^{\mathbf{J}_1(m)}  \bar{r}'_{\mathbf{J}_1(m)m} = 0 \quad \quad \quad \text{for $m\neq m^*$};\\
w^T G^{m^*}= w^T (S^{\eta} \bar{r}'_{\eta m^*} + S^{j_2} \bar{r}'_{j_2 m^*})=0.\\
\end{cases}
\end{equation*}
In particular, this implies
\begin{align}\label{lemnonlem}
\begin{cases}
 w^T S^{\mathbf{J}_1(m)} = 0 \quad \quad \quad\quad\quad\quad \quad \quad &\text{for $m\neq m^*$};\\
 w^T (S^{\eta} r'_{\eta m^*} + S^{j_2} r'_{j_2 m^*})=0 \quad\quad &\text{if and only if } \alpha_{\mathbf{J}_1} r'_{\eta m^*} =\alpha_{\mathbf{J}_2} r'_{j_2m^*}.
\end{cases}
\end{align}
Moreover, \eqref{detresc0} and $G(\bar{\mathbf{r}}') v = 0$ imply 
$$v_{m^*} \neq 0.$$

Let us now fix $\bar{x}$ and $(\mathbf{a}, \mathbf{b})= (\mathpzc{a}, \mathpzc{b})$ as in \eqref{posb1}, \eqref{bconst}, \eqref{a}. For this parameter choice, a straightforward computation yields:
\begin{equation*}
\begin{cases}
f''_{\mathbf{J}_i(m)mm}=-\frac{2r'_{\mathbf{J}_i(m) m}}{x_m} + \frac{(r'_{\mathbf{J}_i(m) m})^2}{r_{\mathbf{J}_i(m)}}\bigg(1+\frac{1}{s^{\mathbf{J}_i(m)}_m}\bigg);\\
f''_{\mathbf{J}_i(m)mn}=f''_{\mathbf{J}_i(m)nm}=\frac{r'_{\mathbf{J}_i(m)m}\; r'_{\mathbf{J}_i(m)n}}{r_j}
\end{cases}
\end{equation*}
For every $\mathbf{J}_i, m, n$, \eqref{mix0} implies that $f''_{\mathbf{J}_i(m)mn}=0$ and hence all mixed second derivatives $g''_{mn}$ are zero, at $\varepsilon=0$. This yields
\begin{equation*}
\begin{split}
w^T \; \partial^2_xg[v,v] &= w^T \sum_{m} \sum_{n} g''_{mn} v_m v_n\\
& = w^T (S^{\eta}f''_{\eta m^* m^*} +S^{j_2}f''_{j_2 m^* m^*})(v_{m^*})^2 + \sum_{m \neq m^*} w S^{\mathbf{J}_1(m)}f''_{\mathbf{J}_1 mm} (v_m)^2\\
&=w^T (S^{\eta}f''_{\eta m^* m^*} +S^{j_2}f''_{j_2 m^* m^*})(v_{m^*})^2.
\end{split}
\end{equation*} 
Via \eqref{lemnonlem}, thus, 
$$w^T \; \partial^2_xg[v,v]\neq 0 \quad\quad\quad \text{if and only if} \quad\quad\quad  \alpha_{\mathbf{J}_1} f''_{\eta m^* m^*}\neq-\alpha_{\mathbf{J}_2} f''_{j_2 m^* m^*}.$$
With \eqref{conddet0} in mind, we compute:
\begin{equation*}
\begin{split}
&\alpha_{\mathbf{J}_1} f''_{\eta m^* m^*}+\alpha_{\mathbf{J}_2} f''_{j_2 m^* m^*} \\
&=  \alpha_{\mathbf{J}_1} \bigg(-\frac{2r'_{\eta m^*}}{x_{m^*}} + \frac{(r'_{\eta m^*})^2}{r_{\eta}}\bigg(1+\frac{1}{s^{\eta}_{m^*}}\bigg) \bigg)+\alpha_{\mathbf{J}_2} \bigg(-\frac{2r'_{j_2 m^*}}{x_{m^*}} + \frac{(r'_{j_2 m^*})^2}{r_{j_2}}\bigg(1+\frac{1}{s^{j_2}_{m^*}}\bigg) \bigg)\\
&=-2\frac{\alpha_{\mathbf{J}_1}r'_{\eta m^*}+\alpha_{\mathbf{J}_2}r'_{j_2 m^*}}{x_{m^*}}+\frac{ \alpha_{\mathbf{J}_1}   (r'_{\eta m^*})^2}{r_\eta}\bigg(1+\frac{1}{s^{\eta}_{m^*}}\bigg) +\frac{\alpha_{\mathbf{J}_2}  (r'_{j_2 m^*})^2}{r_{j_2}}\bigg(1+\frac{1}{s^{j_2}_{m^*}}\bigg)\\
&= - \frac{ \alpha_{\mathbf{J}_2}  r'_{\eta m^*} r'_{j_2 m^*}}{r_\eta}\bigg(1+\frac{1}{s^{\eta}_{m^*}}\bigg) -\frac{\alpha_{\mathbf{J}_1} r'_{\eta m^*} r'_{j_2 m^*}}{r_{j_2}}\bigg(1+\frac{1}{s^{j_2}_{m^*}}\bigg)\\
&=-r'_{\eta m^*} r'_{j_2 m^*} \bigg( \frac{ \alpha_{\mathbf{J}_2} }{r_\eta}\bigg(1+\frac{1}{s^{\eta}_{m^*}}\bigg) +\frac{\alpha_{\mathbf{J}_1}}{r_{j_2}} \bigg(1+\frac{1}{s^{j_2}_{m^*}}\bigg) \bigg) ,
\end{split}
\end{equation*}
which is nonzero if and only if 
$$\frac{\alpha_{\mathbf{J}_2}}{r_\eta}\bigg(1+\frac{1}{s^{\eta}_{m^*}}\bigg) \neq - \frac{\alpha_{\mathbf{J}_1}}{r_{j_2}}\bigg(1+\frac{1}{s^{j_2}_{m^*}}\bigg).$$
By continuity, $w^T \; \partial^2_xg[v,v]\neq 0$ also for small positive $\varepsilon > 0$. Hence, we find Michaelis-Menten functions $f(x, \lambda)$ satisfying all conditions (SN1)--(SN3) of Theorem \ref{SN}.
\endproof

\proof[Proof of Theorem \ref{SNH}]
For Hill kinetics \eqref{MM}, we can use parameters $\mathbf{a},\mathbf{b},\mathbf{c}$. We proceed analogously as in the proof of Theorem \ref{SNMM}, considering the $\varepsilon$-rescaling of the variables $r'_{jm}$ of Lemma \ref{mindistsol}, at the limit $\varepsilon=0$. Again not to overload notation, we omit the explicit dependency on $\varepsilon$. Differently from Theorem \ref{SNMM}, however, we do not have any assumption on the distance $d$.  We again fix $(\mathbf{a}(\mathbf{c}), \mathbf{b}(\mathbf{c}))= (\mathpzc{a}(\mathbf{c}), \mathpzc{b}(\mathbf{c}))$ as in \eqref{posb1}, \eqref{bconst}, \eqref{a}, now explicitly including the dependence on parameters $\mathbf{c}$. We compute again the second derivatives of \eqref{MM} and obtain
\begin{equation*}
\begin{cases}
f''_{jmm}(\mathbf{c})=-\frac{2 c^j_m r'_{j m}}{x_m} + \frac{(r'_{jm})^2}{r_{j}}\bigg(1+\frac{c^j_m x^{c^j_m-1}}{s^{j}_m}\bigg);\\
f''_{jmn}(\mathbf{c})=f''_{jnm}=\frac{r'_{jm}\; r'_{jn}}{r_j}.
\end{cases}
\end{equation*}
At the limit $\varepsilon=0$, $r'_{jm}=0$ implies $f''_{jmm}=0$, $f''_{jmn}=0$, for any choice $\mathbf{c}$ and any $n$. Moreover, a parameter $c^j_m$ appears only in the second derivative $f''_{jmm}$. In particular, the mixed derivatives $f''_{jmn}$ do not depend on the parameters $\mathbf{c}$:
$$\frac{\partial f''_{jmn}(\mathbf{c})}{\partial c^h_k}\equiv 0, \quad\quad\quad\text{for any $h$ and $k$,}$$
while on the contrary, 
$$\frac{\partial f''_{jmm}(\mathbf{c})}{\partial c^h_n}\not\equiv 0,\quad\quad\quad\text{if and only if $h=j$, $n=m$.}$$ 

We focus on the parameter $c^\eta_{m^*}$. At $\varepsilon=0$ we have 
\begin{equation*}
\begin{split}
&w^T \; \partial^2_x g(\bar{x}, \lambda^*)[v,v](c^\eta_{m^*})\\ 
&= w^T \; \sum_{m,n} \frac{\partial^2 g}{\partial x_m \partial x_n}(v_m v_n)(c^\eta_{m^*})\\
&=\sum_{m,n} w^T \; \frac{\partial^2 g}{\partial x_m \partial x_n}(v_m v_n)(c^\eta_{m^*})\\
&=\sum_m w^T \; \frac{\partial^2 g}{(\partial x_m)^2}(v_m)^2(c^\eta_{m^*}) + \sum_{m,n \; m\neq n} w^T \; \frac{\partial^2 g}{\partial x_m \partial x_n}(v_m v_n)\\
&=w^T \; \frac{\partial^2 g}{(\partial x_{m^*})^2}(v_{m^*})^2(c^\eta_{m^*})\\ 
&+ \sum_{m\neq m^*} w^T \; \frac{\partial^2 g}{(\partial x_m)^2}(v_m)^2 + \sum_{m,n \; m\neq n} w^T \; \frac{\partial^2 g}{\partial x_m \partial x_n}(v_m v_n)\\
&=w^T \; S^\eta f''_{\eta m^* m^*}(c^\eta_{m^*})\\ 
&+ w^T \; S^{j_2} f''_{j_2 m^* m^*}  + \sum_{m\neq m^*} w^T \; \frac{\partial^2 g}{(\partial x_m)^2}(v_m)^2 + \sum_{m,n \; m\neq n} w^T \; \frac{\partial^2 g}{\partial x_h \partial x_k}(v_m v_n)
\end{split}
\end{equation*}
We define 
$$K:=w^T \; S^{j_2} f''_{j_2 m^* m^*}  + \sum_{m\neq m^*} w^T \; \frac{\partial^2 g}{(\partial x_m)^2}(v_m)^2 + \sum_{m,n  \;m\neq n} w^T \; \frac{\partial^2 g}{\partial x_m \partial x_n}(v_m v_n).$$  The constant $K$ does not depend on $c^\eta_{m^*}$, and hence $c^\eta_{m^*}$ can be used to nudge $w^T \; \partial^2_xg[v,v]$ away from the degeneracy in the following way. Let us arbitrarily pick a choice $\bar{\mathpzc{c}}$. For $(\mathpzc{a}(\bar{\mathpzc{c}}), \mathpzc{b}(\bar{\mathpzc{c}}), \bar{\mathpzc{c}})$, if the equilibrium $\bar{x}$ possesses a singular Jacobian $G$ with left kernel vector $w$ and right kernel vector $v$ such that
$$w^T \; \partial^2_xg[v,v](\bar{\mathpzc{c}})\neq0,$$
we are done. Otherwise, let us assume 
$$w^T \; \partial^2_xg[v,v](\bar{\mathpzc{c}})= 0.$$
Let us choose $\tilde{\mathpzc{c}}$ such that 
\begin{equation*}
\begin{cases}
\tilde{\mathpzc{c}}^j_m = \bar{\mathpzc{c}}^j_m \quad\quad\quad \text{for $(j, m) \neq (\eta, m^*)$};\\
\tilde{\mathpzc{c}}^j_m \neq \bar{\mathpzc{c}}^j_m \quad\quad\quad \text{for $(j, m) = (\eta, m^*)$}.\\
\end{cases}
\end{equation*}
As previously noted, the parameter $c^\eta_{m^*}$ appears only in $f''_{\eta m^* m^*}$ and hence
\begin{equation*}
\begin{split}
w^T \; \partial^2_x g(\bar{x}, \lambda^*)[v,v] (\tilde{\mathpzc{c}}) - w^T \; \partial^2_x g(\bar{x}, \lambda^*)[v,v] (\bar{\mathpzc{c}})\\
= w^T \; S^\eta f''_{\eta m^* m^*}(\tilde{\mathpzc{c}}^\eta_{m^*})-w^T \; S^\eta f''_{\eta m^* m^*}(\bar{\mathpzc{c}}^\eta_{m^*}) \neq 0,
\end{split}
\end{equation*}
implying $$w^T \; \partial^2_x g(\bar{x}, \lambda^*)[v,v] (\tilde{\mathpzc{c}})\neq 0.$$
By continuity, this extends to small $\varepsilon > 0$.
\endproof

\textbf{Proof of Theorem \ref{SNmy} and Corollary \ref{mindist1}}\\
\nopagebreak

\proof[Proof of Theorem \ref{SNmy}]
The Theorem is just a corollary of Theorem \ref{SNH}. Hill kinetics is indeed a specific example of monotone chemical functions, which proves the theorem.
\endproof

\proof[Proof of Corollary \ref{mindist1}]
We have only to check that the pair $(\tilde{\mathbf{J}}_1,\tilde{\mathbf{J}}_2)$ satisfies Definition \ref{snpairs} of SN-pair of Child Selections of $\tilde{\mathbf{\Gamma}}$. Condition (1) is satisfied by assumption. Condition (2) is inherited from the SN-pair of Child Selections $\mathbf{J}_1, \mathbf{J}_2$ of $\mathbf{\Gamma}$. In fact, note that
$$S^{\tilde{\mathbf{J}_3}}=
\begin{pmatrix}
S^{\mathbf{J}_1} & B \\
0 & D
\end{pmatrix} \quad\quad\quad \text{and} \quad\quad\quad S^{\tilde{\mathbf{J}_4}}=
\begin{pmatrix}
S^{\mathbf{J}_2} & B \\
0 & D
\end{pmatrix},$$ 
where 
$$\begin{pmatrix}
B\\
D
\end{pmatrix}$$
indicates the stoichiometric matrix of the Child Selection $S^{\tilde{\mathbf{J}_3}}$ and $S^{\tilde{\mathbf{J}_4}}$ relative to the reactions $\tilde{\mathbf{J}_3}(m)=\tilde{\mathbf{J}_4}(m)$ for $m \in \tilde{\mathbf{M}}\setminus \mathbf{M}$. Hence,
\begin{equation*}
\begin{cases}
0\neq \alpha_{\tilde{\mathbf{J}}_1} =\operatorname{det}S^{\tilde{\mathbf{J}_3}}=\alpha_{\mathbf{J}_1} \operatorname{det}D\\
0\neq \alpha_{\tilde{\mathbf{J}}_2}=\operatorname{det}S^{\tilde{\mathbf{J}_4}}=\alpha_{\mathbf{J}_2} \operatorname{det}D\\
\end{cases}
\end{equation*}
implying $\alpha_{\tilde{\mathbf{J}}_1}\alpha_{\tilde{\mathbf{J}}_2}<0$.  Condition 3 is satisfied by considering the same species $m^*$ and the same nonzero Partial Child Selection $\mathbf{J}^{\vee m^*}$ for the SN-pair $\mathbf{J}_1, \mathbf{J}_2$ of $\mathbf{\Gamma}$ and extending it to a nonzero Partial Child Selection $\tilde{\mathbf{J}}^{\vee m^*}$ of $\tilde{\mathbf{\Gamma}}$ defined as
\begin{equation*}
\begin{cases}
\tilde{\mathbf{J}}^{\vee m^*}(m)=\mathbf{J}^{\vee m^*}(m) \quad\quad\quad\quad\;\quad\quad \text{for }m\in \mathbf{M}\setminus \{m^*\}\\
\tilde{\mathbf{J}}^{\vee m^*}(m)=\tilde{\mathbf{J}}_1(m)=\tilde{\mathbf{J}}_2(m) \quad\quad\quad\;\; \text{for }m \in \tilde{\mathbf{M}}\setminus \mathbf{M}
\end{cases}.
\end{equation*}
For $\tilde{\mathbf{J}}^{\vee m^*}$ it holds:
$$\beta_{\tilde{\mathbf{J}}^{\vee m^*}} = \operatorname{det} S^{\tilde{\mathbf{J}}^{\vee m^*}}_{\vee m^*}=\operatorname{det} S^{\mathbf{J}^{\vee m^*}}_{\vee m^*} \operatorname{det}D=\beta_{\mathbf{J}^{\vee m^*}}\operatorname{det}D \neq 0,$$
concluding the proof.
\endproof

\section*{Acknowledgments}
I am deeply indebted to Bernold Fiedler for many inspiring discussions. Jia-Yuan Dai helped improve the paper with useful comments.

\bibliographystyle{siamplain}

\end{document}